\documentclass[12pt]{article}

\title{Computation of quantile sets for bivariate ordered data}

\author{
{Andreas H. Hamel\footnote{Free University of Bozen, Faculty of Economics and Management, University Square 1, 39031 Bruneck, Italy, \href{mailto:andreas.hamel@unibz.it}{andreas.hamel@unibz.it}}, Daniel Kostner\footnote{Free University of Bozen, Faculty of Economics and Management, University Square 1, 39100 Bozen, Italy, \href{mailto:daniel.kostner@economics.unibz.it}{daniel.kostner@unibz.it} }}}
\date{}

\usepackage{array} 

\usepackage{verbatim} 
\usepackage{float, graphicx, caption, subcaption}
\usepackage{amsmath, amssymb, enumerate, color}
\usepackage{algpseudocode,algorithm,algorithmicx}

\usepackage{url}

\usepackage{eurosym}

\newtheorem{theorem}{Theorem}
\newtheorem{corollary}[theorem]{Corollary}
\newtheorem{remark}[theorem]{Remark}
\newtheorem{lemma}[theorem]{Lemma}
\newtheorem{definition}[theorem]{Definition}
\newtheorem{proposition}[theorem]{Proposition}
\newtheorem{example}[theorem]{Example}

\numberwithin{equation}{section}  
\numberwithin{figure}{section}    
\numberwithin{table}{section}     
\numberwithin{theorem}{section}

\newcommand{\of}[1]{\ensuremath{\left( #1 \right)}}

\newcommand{\cb}[1]{\ensuremath{ \left\{ #1 \right\} }}
\newcommand{\sqb}[1]{\ensuremath{ \left[ #1 \right] }}

\newcommand{\bs}{\backslash}

\newcommand{\pend}{ \hfill $\square$ \medskip}

\newcommand{\eps}{\ensuremath{\varepsilon}}

\newcommand{\R}{\mathrm{I\negthinspace R}}

\newcommand{\N}{\mathrm{I\negthinspace N}}

\newcommand{\co}{{\rm co \,}}

\newcommand{\Int}{{\rm int\,}}

\usepackage[
colorlinks=true, linkcolor=blue, citecolor=blue, urlcolor=blue, filecolor=blue
auskommentieren und unten pdfborder aktivieren]{hyperref}

\definecolor{color0}{gray}{.50}
\definecolor{color1}{rgb}{0,.2,.8}
\definecolor{color2}{rgb}{1,.2,0}
\definecolor{color3}{rgb}{.8,.5,1}
\newcommand{\f}{\color{color1}}
\newcommand{\ff}{\color{color2}}

\begin{document}
\maketitle

\begin{abstract}
Algorithms are proposed for the computation of set-valued quantiles and the values of the lower cone distribution function for bivariate data sets. These new objects make data analysis possible involving an order relation for the data points in form of a vector order in two dimensions. The bivariate case deserves special attention since two-dimensional vector orders are much simpler to handle than such orders in higher dimensions. Several examples illustrate how the algorithms work and what kind of conclusions can be drawn with the proposed approach.
\end{abstract}



\section{Introduction}

Algorithms are presented for computing the quantile sets for bivariate random variables as well as the values of the corresponding lower cone distribution function in the presence of an order relation for their values. Such quantiles in the general multivariate case have been defined in \cite{HamelKostner18JMVA} as a common generalization of univariate quantiles and Tukey's (halfspace) depth regions. Likewise, the lower cone distribution function is a common generalization of the univariate cumulative distribution function and Tukey's (halfspace) depth function. It can also be used as a ranking function for multi-criteria decision making problems \cite{Kostner20MMOR}.

Moreover, it was shown in \cite{AraratHamel20TPA} that set-valued quantiles and the lower cone distribution functions form Galois connections between complete lattices of sets and the interval [0,1] of real number. This generalizes a property which is straightforward and well-known in the univariate case, but has never been discussed with respect to depth functions and depth regions. It is also shown in \cite{AraratHamel20TPA} that set-valued quantiles characterize the distribution of a random set extension of the original random variable as well as its capacity functional.

The bivariate case deserves special attention since convex cones in $\R^2$ have a very simple representation (every closed convex cone is polyhedral, i.e., the intersection of a finite number of halfspaces---this number being 1 or 2 in almost all cases) and, using this, the computations can be done much faster than in the general $\R^d$-valued case: there are polyhedral cones with arbitrary many facets already in $\R^3$.

Algorithms for the bivariate location depth and corresponding depth regions were given in \cite{RousseeuwRuts96JRSS, RutsRousseeuw96CSDA}. Algorithms for depth functions and regions in general dimensions can be found, for example, in \cite{DyckerhoffMozharovskyi16CSDA, Liu17CS, LiuMoslerMozharovskyi19JCGS, LiuZuo14CSSC, PaindaveineSiman12CSDA}. These references are mainly concerned with Tukey depth functions/regions, i.e., they do not take an order relation for the value of the random variable into account. The reader may compare \cite{BelloniWinkler11AS} which is one of the very few references dealing with statistics for multivariate ordered data.

On the other hand, an order relation is often present and intuitive since decision makers have preferences or the impact of some events is clearly preferred over the ones of others. A few examples illustrating this feature are discussed below: hurricane scales, hail insurance and human resource management. It is beyond the scope of this paper to give an exhaustive statistical analysis of these example; they are used as showcases for the type of conclusions which can be drawn with the approach initiated in \cite{HamelKostner18JMVA}, in particular what sets it apart from a mere depth function approach. 

The paper is organized as follows. In the next section, vector orders in $\R^2$ are reviewed. Section \ref{SecEmpCDF-Quan} includes the definition of the main concepts and preparatory results. The algorithms are presented in Section \ref{SecBivariate} while in Section \ref{SecEx} several examples are discussed including a new view on hurricane scales and the problem of finding best candidates for tasks/jobs.

\section{Vector preorders in two dimensions}
\label{SubSecPreorders2D}

The basic assumption is that there is a preference relation for the two-dimensional data points in form of a vector preorder, i.e., a reflexive and transitive relation which is compatible with the algebraic operations in $\R^2$. Such vector preorders are in one-two-one correspondence with convex cones $C \subseteq \R^2$ including $0 \in \R^2$ via
\begin{equation}
\label{EqConePreorder}
y \leq_C z \quad \Leftrightarrow \quad z - y \in C
\end{equation}
(see, for example, \cite[Chap. 8]{AliprantisBorder06Book3rd}). A convex cone $C \subseteq \R^2$ is a set satisfying $s C \subseteq C$ for all $s > 0$ and $C + C \subseteq C$.

In the following, it is assumed that the cone $C$ generating the preorder via \eqref{EqConePreorder} is closed. Such vector preorders in $\R^2$ have some special features compared to the case $d >2$. Only the following cases are possible:

1. The cone is a linear subspace which is either $C =\{0\}$, or $C = \R^2$, or a straight line $C = L$ in $\R^2$.

2. The cone is a ray: $C = \cb{s b \mid s \geq 0}$ for some $b \in \R^2\bs\{0\}$.

3. The cone is generated by two linearly independent vectors: $C = \cb{s_1b^1 + s_2b^2 \mid s_1, s_2 \geq 0}$ for some $b^1, b^2 \in \R^2$ which are linearly independent.

4. The cone is a closed (homogeneous) halfspace:  $C = \cb{s_1b - s_2b + s_3\bar z \mid s_1, s_2, s_3 \geq 0}$ for two linearly independent vectors $b, \bar z \in \R^2\bs\{0\}$.

There are interesting non-closed cones such as the lexicographic ordering cone even in $\R^2$; such cases require a different type of analysis (since the bipolar theorem does not apply) and therefore, they will not be considered here.

The case $C=\{0\}$ leads to the Tukey halfspace depth function and regions; it is dealt with, e.g., already in \cite{RutsRousseeuw96CSDA}. While the case $C=\R^2$ is trivial, the case $C = L$ will not be discussed in this paper. Finally, if $C$ is a closed halfspace, there is $w \in \R^2\bs\{0\}$ such that $C = H^+(w) = \{z \in \R^2 \mid w^Tz \geq 0\}$; in this case, the order $\leq_{H^+(w)}$ is a total preorder (a reflexive, transitive relation such that either $y \leq_{H^+(w)} z$ or $z \leq_{H^+(w)} y$ or both) and the situation can be reduced to the univariate case. 

In this paper, the main subject is the case of a closed convex pointed cone with non-empty interior, i.e., \#3 above and case \#4 will appear as an intermediate step.

If $C$ is generated by two linearly independent vectors $b^1, b^2 \in \R^2$, then $C$ is the intersection of exactly two halfspaces, i.e., there are $v^1, v^2 \in \R^2\bs\{0\}$ linearly independent such that 
\[
C = H^+(v^1) \cap H^+(v^2)
\]
(choose $v^1$ and $v^2$ orthogonal to  $b^1$ and $b^2$, respectively, such that $(v^1)^Tb^2 \geq 0$, $(v^2)^Tb^1 \geq 0$).

It is assumed in the following that $b^1, b^2$ as well as $v^1, v^2$ are known, and the two sets $\cb{b^1, b^2}$ and $\cb{v^1, v^2}$ are called a V-representation and an H-representation of $C$, respectively. The set 
\[
C^+ = \cb{w \in \R^2 \mid \forall z \in C \colon w^Tz \geq 0} = \cb{w \in \R^2 \mid w^Tb^1, w^Tb^2 \geq 0}
\]
is called the (positive) dual cone of $C$ (always a closed convex cone). Under the given assumptions,
\[
C^+ = \bigcup_{t \geq 0}tB^+,
\]
where $B^+ = \cb{sv^1 + (1-s)v^2 \mid s \in [0,1]}$ is a base of $C^+$, i.e., for each $w \in C^+\bs\{0\}$ there are unique $t > 0$ and $b \in B^+$ such that $w = tv$.

\section{Empirical cone distribution functions and quantiles}
\label{SecEmpCDF-Quan}

In this section, we give the definitions of lower cone distribution functions and associated quantiles for bivariate random variables in case of a finite data sets. Let $\tilde X = \cb{x^1, x^2, \ldots, x^N} \subseteq \R^2$ be a finite collection of data points which could be a sample of a random variable. The following definition provides the bivariate empirical counterpart to the concepts from \cite{HamelKostner18JMVA}. Compare also \cite{Kostner20MMOR} for the $\R^d$-valued case with applications to a multi-criteria decision making problem.

\begin{definition}
\label{DefEmpCDF} The functions $F_{X, w} \colon \R^2 \to [0,1]$ for $w \in B^+$ and $F_{X, C} \colon \R^2 \to [0,1]$ defined by
\begin{align}
F_{\tilde X, w}(z) & = \frac{1}{N}\#\cb{x \in \tilde X \mid x \in z - H^+(w)} \quad \text{and} \label{EqW-LDF}\\
F_{\tilde X, C}(z) & = \min_{w \in B^+} F_{\tilde X, w}(z) = \frac{1}{N}\min_{w \in B^+} \#\cb{x \in \tilde X \mid x \in z - H^+(w)}. \label{EqC-LDF}
\end{align}
are called empirical lower $w$-distribution function and empirical lower $C$-distribution function, respectively, for the data set $\tilde X$.
\end{definition}

The functions $w\text{-}depth(z; \tilde X) := N \cdot F_{\tilde X, w}(z)$ and $c\text{-}depth(z; \tilde X) := N \cdot F_{\tilde X, C}(z)$ are called the {\em $w$-location depth} and the {\em cone location depth} for $\tilde X$. 

The functions $w\text{-}depth$ and $c\text{-}depth$ can be interpreted as follows. For each point $z \in \R^2$, the $w$-location depth gives the number of data points which are dominated by $z$ with respect to the total preorder generated by $H^+(w)$, i.e., data points $x \in \tilde X$ satisfying $w^\top x \leq w^\top z$. The cone location depth of $z \in \R^2$ gives the minimal number of data points which are dominated by $z$ with respect to all total preorders generated by $H^+(w)$ for $w \in C^+$. Thus, a  point $z \in \R^2$ dominates at least $c\text{-}depth(z; \tilde X)$ data points with respect to the total preorder generated by $H^+(w)$ for all $w \in C^+$, i.e., no matter which weighted average with weights from $C^+$ is taken. Data points which are higher ranked than others are ``deeper" in the sense that they improve with respect to more $w$'s at the same time.

\begin{proposition}
\label{PropMono}
(1) $w^\top y \leq w^\top z$ implies $w\text{-}depth(y; \tilde X) \leq w\text{-}depth(z; \tilde X)$;

(2) $y \leq_C z$ implies $c\text{-}depth(y; \tilde X) \leq c\text{-}depth(z; \tilde X)$.
\end{proposition}

{\sc Proof.} This follows directly from the monotonicity property of $F_{\tilde X, w}$ and $F_{\tilde X, C}$, respectively, in \cite[Proposition 1 (b)]{HamelKostner18JMVA}. \pend

\medskip
The following example shows that the cone location depth (as well as $F_{\tilde X, C}$) can be understood as a ranking function for the data points which reflects the order $\leq_C$. This seems to be very much in the spirit of Tukey's original work. This example also shows that points which are non-comparable with respect to the order $\leq_C$ can have the same or very different cone location depths.

\begin{figure}[H]
\centering
\includegraphics[width=1\textwidth]{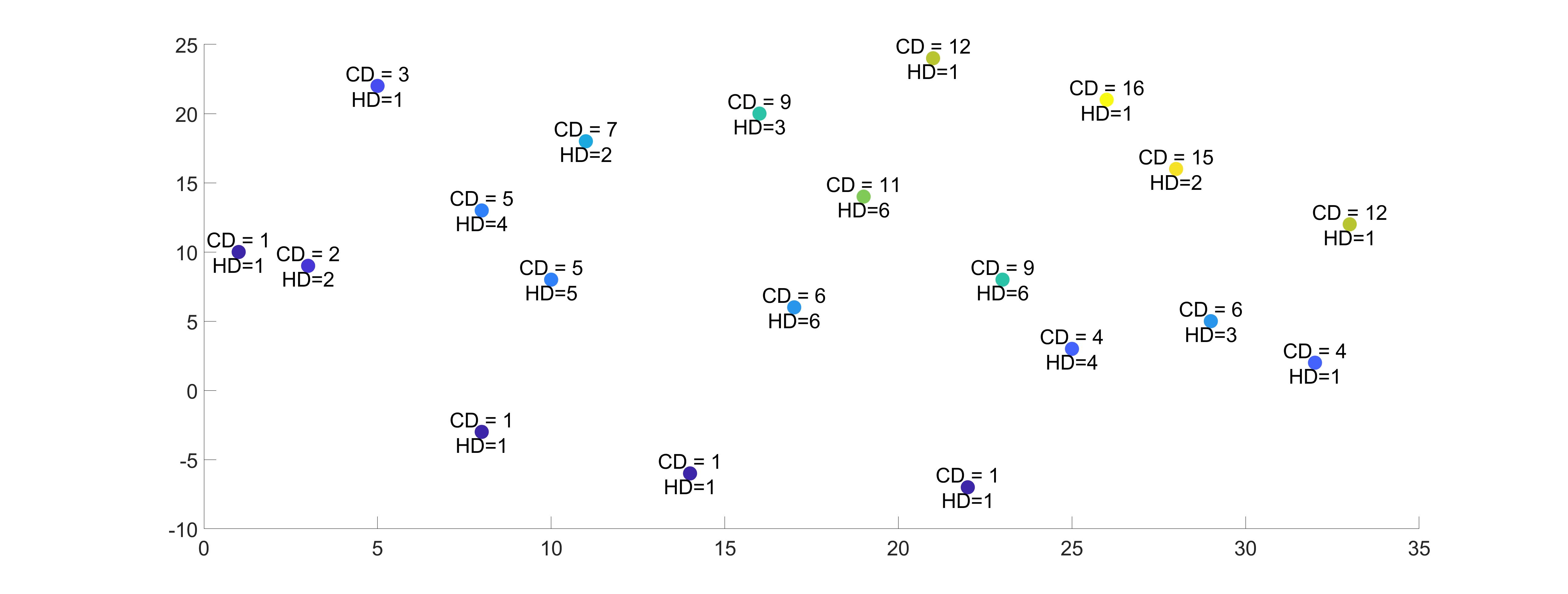}
\caption{$C = \R^2_+, N = 20, p = 0.2, \lceil Np \rceil = 4$}
\label{fig:DepthPoints}
\end{figure}

\begin{example}
\label{ExDepth}
For every data point in Figure \ref{fig:DepthPoints} Tukey's (halfspace) depth HD and the cone location depth CD are computed. One may already realize that the cone location depth ``follows the cone" (increases in directions in which the cone ``opens") whereas the halfspace depth increases toward the center of the data cloud. This means that the cone location depth ranks the data points taking into account the order generated by the cone. This is a new feature not captured by depth functions.
\end{example}

Empirical quantiles are set-valued functions, i.e., they map into the power set $\mathcal P(\R^2)$, the set of all subsets of $\R^2$ including the empty set.

\begin{definition}
\label{DefEmpQuant} 
The empirical $w$-quantile function $Q^-_{\tilde X, w} \colon [0,1] \to \mathcal P(\R^2)$ and the empirical $C$-quantile function $Q^-_{\tilde X, C}\colon [0,1] \to \mathcal P(\R^2)$ associated to $\tilde X$ and $C$ are defined by
\[
Q^-_{\tilde X, w}(p )  = \cb{z \in \R^2 \mid  F_{\tilde X, w}(z) \geq p} \quad \text{and} \quad
Q^-_{\tilde X, C}(p )  = \cb{z \in \R^2 \mid F_{\tilde X, C}(z) \geq p},
\]
respectively.
\end{definition}

The definitions of $F_{\tilde X, w}$ and $F_{\tilde X, C}$ immediately yield
\begin{align}
Q^-_{\tilde X, w}(p )  & = \cb{z \in \R^2 \mid \#\{x \in \tilde X \mid x \in z - H^+(w)\} \geq \lceil Np \rceil} 
	 \label{EqDiscreteW-Q}  \\
Q^-_{\tilde X, C}(p )  & = \cb{z \in \R^2 \mid \min_{w \in B^+}\#\{x \in \tilde X \mid x \in z - H^+(w)\} \geq \lceil Np \rceil} \label{EqDiscreteC-Q}
\end{align}
where $\lceil Np \rceil$ is the value of the ceiling function at $Np$ defined by $\lceil r \rceil = \min\cb{k \in \N \mid r \leq k}$ (the least natural number which is greater than or equal to $r \in \R$). Clearly, 
\begin{align*}
Q^-_{\tilde X, w}(p ) & = \cb{z \in \R^d \mid w\text{-}depth(z; \tilde X)  \geq \lceil Np \rceil} \quad \text{and} \\ 
Q^-_{\tilde X, C}(p ) & =   \cb{z \in \R^d \mid c\text{-}depth(z; \tilde X)  \geq \lceil Np \rceil}.
\end{align*}

Let $\tilde Y \subseteq \R$ be a finite univariate data set. We denote its empirical lower quantile by $q^-_{\tilde Y}(p ) = \min\{\bar y \in \tilde Y \mid \#\{y \in \tilde Y \mid y \leq \bar y\} \geq \lceil Np \rceil\}$. With this notation, one has
\begin{align}
\forall p \in (0,1) \colon Q^-_{\tilde X, w}(p ) & = \cb{z \in \R^2 \mid w^\top z \geq  q^-_{w^\top\tilde X}(p )} \quad \text{and} \quad
\label{EqWScalarization}\\
\forall p \in (0,1) \colon Q^-_{\tilde X, C}(p ) & = \bigcap_{w \in B^+}\cb{z \in \R^d \mid w^\top z \geq  q^-_{w^\top\tilde X}(p )} \label{EqC-Scalarization}
\end{align}
as in the general (multi-dimensional) case (see \cite[Proposition 6]{HamelKostner18JMVA}).

Next, we formally state that at least one data point must be on the boundary of each $w$-quantile.

\begin{proposition}
\label{PropDataPtBoundary}
If $p \in [0,1]$, $w \in B^+$ and $Q^-_{\tilde X, w}(p ) \not\in \{\R^2, \emptyset\}$, then there is $x(w, p) \in \tilde X$ such that 
\begin{equation}
\label{EqW-Quantile}
Q^-_{\tilde X, w}(p ) = x(w, p) + H^+(w).
\end{equation}
Moreover, the following three conditions are equivalent for $\tilde x \in \tilde X$:

(a) $Q^-_{\tilde X, w}(p ) = \tilde x + H^+(w)$.

(b) One has
\begin{align}
& \#\{x \in \tilde X \mid x \in \tilde x - H^+(w)\} \geq \lceil Np \rceil, \label{EqNumNegWQuantile} \\
& \#\{x \in \tilde X \mid x \in \tilde x - \Int H^+(w)\} < \lceil Np \rceil. \label{EqNumIntNegWQuantile}
\end{align}

(c) $q^-_{w^\top \tilde X}(p ) = w^\top \tilde x$.
\end{proposition}

{\sc Proof.} Since $Q^-_{\tilde X, w}(p )$ is a closed halfspace with normal $w$, there is a point $y \in \R^d$ such that $Q^-_{\tilde X, w}(p ) = y + H^+(w)$. By \eqref{EqDiscreteW-Q}, $\#\cb{x \in \tilde X \mid x \in y - H^+(w)} \geq \lceil Np\rceil$. If there would be no data point on the boundary of this halfspace, one even had $\#\cb{x \in \tilde X \mid x \in y - \Int H^+(w)} \geq \lceil Np\rceil$. Then, there would exist $y' \in y - \Int H^+(w)$ with $\#\cb{x \in \tilde X \mid x \in y' - \Int H^+(w)} \geq \lceil Np\rceil$, hence $y' \in Q^-_{\tilde X, w}(p )$, but $y' \not\in y + H^+(w)$, a contradiction.

(a) $\Rightarrow$ (b): If (a) is true, then $\tilde x \in Q^-_{\tilde X, w}(p )$ and \eqref{EqNumNegWQuantile} follows from \eqref{EqDiscreteW-Q}. If \eqref{EqNumIntNegWQuantile} would not be true, then, by a similar argument as in the first part of the proof, $\tilde x$ would not be in $Q^-_{\tilde X, w}(p )$.

(b) $\Rightarrow$ (a): Assume $z \in \tilde x + H^+(w)$. Then $z - H^+(w) \supseteq \tilde x - H^+(w)$, hence
\[
 \#\{x \in \tilde X \mid x \in z - H^+(w)\} \geq  \#\{x \in \tilde X \mid x \in \tilde x - H^+(w)\} \geq \lceil Np \rceil
\] 
which means $z \in Q^-_{\tilde X, w}(p )$. Hence $\tilde x + H^+(w) \subseteq Q^-_{\tilde X, w}(p )$. Conversely, if $z \in Q^-_{\tilde X, w}(p )$, then $\#\{x \in \tilde X \mid x \in z - H^+(w)\}  \geq \lceil Np \rceil$ by \eqref{EqNumNegWQuantile}. This and \eqref{EqNumIntNegWQuantile} imply $z \not\in \tilde x - \Int H^+(w)$ (otherwise, $z - H^+(w) \subseteq \tilde x - \Int H^+(w)$ which leads to a contradiction). But then $z \in \tilde x + H^+(w)$, hence $Q^-_{\tilde X, w}(p ) \subseteq \tilde x + H^+(w)$.

(a) $\Leftrightarrow$ (c) Both directions follow from \eqref{EqWScalarization}. \pend

\medskip For the following result, some notation is needed which is also used in the remainder of the paper. For $w \in B^+$ we define the following sets
\begin{align*}
\tilde X^=(w, p) & =  \cb{x \in \tilde X \mid w^\top x = w^\top x(w, p)} \\
\tilde X^\leq(w, p) & =  \cb{x \in \tilde X \mid w^\top x \leq w^\top x(w, p)}. 
\end{align*}
The set $\tilde X^=(w, p)$ includes all data points on the boundary of the shifted halfspace $Q^-_{\tilde X, w}(p ) = x(w, p) + H^+(w)$, whereas $\tilde X^\leq(w, p)$ includes the data points in $x(w, p) - H^+(w)$ with $x(w, p) \in \tilde X$ from \eqref{EqW-Quantile} in both cases.

\medskip {\bf Standing assumption.} In the remainder of the paper, it is assumed that $C^+$ is generated by the two linearly independent vectors $v^1, v^2 \in C^+$, i.e., $\{v^1, v^2\}$ is a V-representation of $C^+$. In this case, the set 
\[
B^+ = \cb{(1-s)v^1 + s v^2 \mid s \in [0,1]}
\]
is a base of $C^+$. Equivalently, $C$ is generated by two linearly independent vectors $b^1, b^2 \in C$, i.e., $\{b^1, b^2\}$ is a V-representation of $C$. 

\begin{proposition}
\label{PropIntNormals}
If $p \in (0,1]$ and $Q^-_{\tilde X, C}(p ) \not\in \{\R^2, \emptyset\}$, then $Q^-_{\tilde X, C}(p )$ can be represented as
\[
Q^-_{\tilde X, C}(p ) = \bigcap_{w \in W(p )} Q^-_{\tilde X, w}(p )
\]
such that $W(p ) \subseteq B^+$ and 
\[
\#\tilde X^=(w, p) \geq 2 \quad\text{and}\quad \# \tilde X^\leq(w, p) \geq \lceil Np\rceil +1
\]
for all $w \in W(p )\bs\{v^1, v^2\}$. In particular, $W(p )$ is a finite set.
\end{proposition}

{\sc Proof.} First, take $w \in B^+\bs\{v^1, v^2\}$. By Proposition \ref{PropDataPtBoundary}, $\#\tilde X^=(w, p) \geq 1$ where
\[
\tilde X^=(w, p) = \cb{x \in \tilde X \mid  w^\top x = w^\top x(w, p)} = \cb{x \in \tilde X \mid Q^-_{\tilde X, w}(p ) = x + H^+(w)}.
\]
Since $w \not\in \{v^1, v^2\}$, there is $s \in (0, 1)$ such that $w = (1-s)v^1 + s v^2$. Assume that $\#\tilde X^=(w, p) = 1$, i.e., $\tilde X^=(w, p)  = \{x(w,p)\}$. Then, there are $\eps_1, \eps_2 > 0$ such that $s_1 := s + \eps_1 \in (0,1)$, $s_2 := s - \eps_2 \in (0,1)$ and the two halfspaces
\[
x(w,p) + H^+(w(s_1)) \quad \text{and} \quad x(w,p) + H^+(w(s_2))
\]
contain exactly the same set of data points as $Q^-_{\tilde X, w}(p ) = x(w,p) + H^+(w)$ where $w(s_i) = (1-s_i)v^1 + s_iv^2 \in B^+$, $i=1,2$. Hence $Q^-_{\tilde X, w(s_i)}(p ) = x(w,p) + H^+(w(s_i))$, $i=1,2$, and
\[
x(w,p) + H^+(w) \supsetneq \of{x(w,p) + H^+(w(s_1))} \bigcap \of{x(w,p) + H^+(w(s_2))}.
\]
This means that $Q^-_{\tilde X, w}(p )$ does not contribute to the intersection in
\begin{equation}
\label{EqQCIntersection}
Q^-_{\tilde X, C}(p ) = \bigcap_{w \in B^+} Q^-_{\tilde X, w}(p ),
\end{equation}
and it is enough to run it over those $w \in B^+$ with $\#\tilde X^=(w, p) \geq 2$ and $v^1, v^2$. 

Secondly, take such a $w \in B^+ \bs\cb{v^1, v^2}$ and assume $\#\tilde X^\leq(w, p) \leq \lceil Np\rceil$. Now, \eqref{EqNumNegWQuantile} implies "=" in this inequality. Pick $\tilde x^1, \tilde x^2 \in \tilde X^=(w, p)$ such that 
\begin{eqnarray}
\label{EqV1Max}
(v^1)^\top \tilde x^1  & = \max\cb{(v^1)^\top x \mid x \in  \tilde X^=(w, p)}, \\
\label{EqV2Max}
(v^2)^\top \tilde x^2  & = \max\cb{(v^2)^\top x \mid x \in  \tilde X^=(w, p)}.
\end{eqnarray}
One has $\tilde x^1 \neq \tilde x^2$ and $(v^1)^\top\tilde x^1 > (v^1)^\top \tilde x^2$ and $(v^2)^\top\tilde x^2 > (v^2)^\top \tilde x^1$ by \eqref{EqV1Max}, \eqref{EqV2Max}, hence
\begin{equation}
\label{EqVMax}
\of{v^1 - v^2}^\top\of{\tilde x^1 - \tilde x^2} > 0.
\end{equation}
Since the data set $\tilde X$ is finite, it is always possible to find $\eps_1, \eps_2 > 0$ and $w^1 := w(\eps_1), w^2 := w(\eps_2)$ such that 
\begin{eqnarray*}
w^1(\eps) & = w + \eps\of{v^1 - v^2} \in B^+ \; \text{for} \; \eps \in [0, \eps_1] \\
w^2(\eps) & = w + \eps\of{v^2 - v^1} \in B^+ \; \text{for} \; \eps \in [0, \eps_2]
\end{eqnarray*}
and the following conditions are satisfied: $\tilde x^i$ is kept on the boundary of the halfspace $\tilde x^i + H^+(w^i(\eps))$ for $i=1,2$, and one has $\#\tilde X^=(w^i(\eps), p) = 1$,  $Q^-_{\tilde X, w^i(\eps)}(p ) = \tilde x^i + H^+(w^i(\eps))$ and $\#\tilde X^\leq(w^i(\eps), p) = \lceil Np\rceil$ for $\eps \in [0, \eps_i)$, and for $w^i$ (i.e., for $\eps = \eps_i$), $i = 1,2$, one has $Q^-_{\tilde X, w^i}(p ) = \tilde x^i + H^+(w^i)$, $\#\tilde X^\leq(w^i, p) \geq \lceil Np\rceil$ and either $w^i = v^i$ or $\#\tilde X^=(w^i, p) \geq 2$. 

The underlying geometrical idea is to turn $w$ in direction $v^1$ and $v^2$, respectively, around $\tilde x^1$ and $\tilde x^2$ until the next data point is hit. The data points in $\tilde X^\leq(w^i(\eps), p)$ with $\eps \in [0, \eps_i)$, $i=1,2$, are exactly the same as in $\tilde X^\leq(w, p)$.

Then, one has
\begin{equation}
\label{EqRedundant}
\tilde x^1 + H^+(w) = \tilde x^2 + H^+(w) \supsetneq \sqb{\tilde x^1 + H^+(w^1)} \cap \sqb{\tilde x^2 + H^+(w^2)}
\end{equation}
which means that $w$ is indeed redundant in the intersection \eqref{EqQCIntersection}. To see this, observe that there is $\hat z \in 
\sqb{\tilde x^1 + H^+(w^1)} \cap \sqb{\tilde x^2 + H^+(w^2)}$ satisfying
\[
\tilde x^1 + H^+(w^1)  = \hat z + H^+(w^1) \quad \text{and} \quad \tilde x^2 + H^+(w^2)  = \hat z + H^+(w^2)
\]
($\hat z$ is the intersection point of the two boundary lines of $\tilde x^1 + H^+(w^1)$, $\tilde x^2 + H^+(w^2)$ and it does not have to be a data point, of course).

{\bf Claim.} $\hat z \in \tilde x^i + \Int H^+(w)$, i.e., $w^T(\hat  z - \tilde x^i) >0$, for $i=1,2$.

Indeed, one has
\[
(w^1)^\top(\hat z - \tilde x^1) = 0, \quad (w^2)^\top(\hat z - \tilde x^2) = 0,
\]
i.e.,
\begin{eqnarray}
\label{EqZetHat1}
\of{w + \eps_1\of{v^1 - v^2}}^\top(\hat z - \tilde x^1) = 0 \\
\label{EqZetHat2}
\of{w - \eps_2\of{v^1 - v^2}}^\top(\hat z - \tilde x^2) = 0.
\end{eqnarray}
Multiplying the second equation by $-1$ and adding the result to the first gives
\[
\eps_1\of{v^1 - v^2}^\top(\hat z - \tilde x^1) + \eps_2\of{v^1 - v^2}^\top(\hat z - \tilde x^2) = 0,
\]
so one of the two parts of the sum must be $\leq 0$, the other $\geq 0$. With the help of \eqref{EqV1Max}, \eqref{EqV2Max} one gets
\begin{align*}
\of{v^1 - v^2}^\top(\hat z - \tilde x^1) & = \of{v^1 - v^2}^\top\hat z - (v^1)^\top\tilde x^1 + (v^2)^\top\tilde x^1 \\
	& < \of{v^1 - v^2}^\top\hat z - (v^1)^\top\tilde x^2 + (v^2)^\top\tilde x^2 \\
	& = \of{v^1 - v^2}^\top(\hat z - \tilde x^2),
\end{align*}
hence
\[
\of{v^1 - v^2}^\top(\hat z - \tilde x^1) < 0, \quad \of{v^1 - v^2}^\top(\hat z - \tilde x^2) > 0.
\]
Together with \eqref{EqZetHat1}, \eqref{EqZetHat2}, this gives
\begin{eqnarray*}
w^\top(\hat z - \tilde x^1) & = -\eps_1\of{v^1 - v^2}^\top(\hat z - \tilde x^1) > 0 \\
w^\top(\hat z - \tilde x^2) & = \eps_2\of{v^1 - v^2}^\top(\hat z - \tilde x^2) > 0,
\end{eqnarray*}
so $\hat z \in \tilde x^i + \Int H^+(w)$ for $i = 1,2$ which proves the claim.

Finally, take
\[
z \in \sqb{\tilde x^1 + H^+(w^1)} \cap \sqb{\tilde x^2 + H^+(w^2)} = \sqb{\hat z + H^+(w^1)} \cap \sqb{\hat z + H^+(w^2)},
\]
i.e.,
\[
(w^1)^\top(z - \hat z) \geq 0, \quad (w^2)^\top(z - \hat z) \geq 0.
\]
The definitions of $w^1$ and $w^2$ yield
\[
\of{w + \eps_1\of{v^1 - v^2}}^\top(z - \hat z) \geq 0, \quad \of{w - \eps_2\of{v^1 - v^2}}^\top(z - \hat z) \geq 0,
\]
hence
\[
w^\top(z - \hat z) \geq \max\{-\eps_1\of{v^1 - v^2}^\top(z - \hat z), \eps_2\of{v^1 - v^2}^\top(z - \hat z)\} \geq 0,
\]
so $z \in \hat z + H^+(w) \subseteq \tilde x^i + \Int H^+(w)$, $i=1,2$, where the inclusion is the claim above. This proves \eqref{EqRedundant} which completes the proof of the proposition. \pend

\begin{remark}
For $i = 1,2$, one can have
\[
\#\tilde X^\leq(w^i, p) = \lceil Np\rceil \quad \text{or} \quad \#\tilde X^\leq(w^i, p) \geq \lceil Np\rceil + 1.
\]
In the first case, $w^i$ is also redundant by Proposition \ref{PropIntNormals}. This is exploited in the algorithm below. In the second, $w^i$ cannot be ruled out by the proposition, but it could still be redundant for the intersection in \eqref{EqQCIntersection}. 
\end{remark}

\begin{remark}
Under the standing assumption, the situation can be reduced to the case when the cone is $C = \R^2_+$. Let $C$ be generated by the two linearly independent vectors $b^1, b^2 \in \R^2$, i.e., $C = \{s_1 b^1 + s_2 b^2 \mid s_1 \geq 0, \, s_2 \geq 0\}$. Then, there is an invertible matrix $A \in \R^{2 \times 2}$ such that $AC = \R^2_+$ and one can use the affine invariance of the cone distribution function $F_{X,C}$ (see \cite[Proposition 2.7]{HamelKostner18JMVA}) to get
\[
\forall z \in  \R^2 \colon F_{X,C}(z) = F_{AX, \R^2_+}(Az).
\] 
Indeed, since $AC = \R^2_+$ if, and only if, $Ab^1 = e^1$ and $Ab^2 = e^2$, one may easily see that
\[
A^{-1} = 
\of{
	\begin{array}{cc}
	b^1_1 & b^2_1 \\ 
	b^1_2 & b^2_2 \\ 
	\end{array}
}
\quad \text{and} \quad
A = \frac{1}{b^1_1b^2_2 - b^1_2b^2_1}
\of{
	\begin{array}{cc}
  b^2_2 & -b^2_1 \\ 
  -b^1_2 & b^1_1 \\ 
\end{array}
}
\]
do the job. Moreover,
\[
w \in C^+ \quad \Leftrightarrow \quad \of{A^{-1}}^\top w \in \R^2_+,
\]
so $\of{A^{-1}}^\top C^+  = \R^2_+$. The procedure now is: first, transform the data and the cone $C$ by $A$; secondly find $F_{AX, \R^2_+}$ and $AQ^-_{X,C}$; finally transform back. Clearly, this idea is restricted to the bivariate case.
\end{remark}

\section{The bivariate algorithms}
\label{SecBivariate}

In this section, we provide some more theoretical background for algorithms which produce the empirical quantiles $Q^-_{\tilde X, w}(p )$ with empirical halfspace depth regions as special cases and the values of $c\text{-}depth$ with the halfspace depth function as a special case, respectively. Pseudocodes of the algorithms will also be given.

\subsection{The rotation step}
\label{SubSecRotation}


The algorithm below is designed such that it starts with $w^1 = v^1$ and then runs through $B^+$ until it hits $v^2$. At an intermediate step, a $w^n \in B^+$ is generated.  There are three cases:

(1) $w^n = v^1$ and $\#X^=(v^1, p) = 1$.

(2) If $w^n \neq v^2$ and $\#X^=(w^n, p) \geq 2$. 

(3) $w^n = v^2$. 
\\
Case (3) serves as a stopping criterion. In case (1) and (2), a permutation of the data points is generated which in turn is used to generate a new $w^{n+1} \in B^+$. Set $K := \lceil Np\rceil$.

\medskip
{\bf Case (1).} The input is $w:= w^1 = v^1 \in B^+$, $x(w, p)$ (the only element in $\#X^=(w, p)$).  Find the permutation $\pi_w$ of $\{1, \ldots, N\}$ such that
\begin{equation}
\label{EqRotPermutation0}
w^\top x^{\pi_w(1)} \leq \ldots < w^\top x^{\pi_w(K)} < \ldots \leq w^\top x^{\pi_w(N)}.
\end{equation}
One has $x^{\pi_w(K)} = x(w, p)$. 

\medskip
{\bf Case (2).} The input is $w:= w^n \in B^+$. Find the permutation $\pi_w$ of $\{1, \ldots, N\}$ such that
\begin{equation}
\label{EqRotPermutation1}
w^\top x^{\pi_w(1)} \leq \ldots \leq w^\top x^{\pi_w(K)} \leq \ldots \leq w^\top x^{\pi_w(N)}
\end{equation}
Find the set
\[
X^=(w, p) = \cb{x \in \tilde X \mid w^\top x = w^\top x^{\pi_w(K)}}
\]
and set $L := \#X^=(w, p)$. Define a new permutation $\pi_{v^2}$ of the points in $X^=(w, p)$ by
\begin{equation}
\label{EqPermEqual}
(v^2)^\top x^{\pi_{v^2}(k_1)} \leq \ldots \leq (v^2)^\top x^{\pi_{v^2}(k_L)}
\end{equation}
of the $L$ data points in $\#X^=(w, p)$.

\begin{proposition}
\label{PropStrictPermut} 
All inequalities in \eqref{EqPermEqual} are strict.
\end{proposition}

{\sc Proof.} Assume to the contrary that $(v^2)^\top x = (v^2)^\top y$ with $x, y \in X^=(w, p)$. One has $w = (1-s)v^1 + sv^2$ for some $s \in [0,1)$. 

First, assume that $s = 0$. Then $w= v^1$ and $(v^1)^\top x= (v^1)^\top y$ from \eqref{EqRotPermutation1} as well as $(v^2)^\top x= (v^2)^\top y$. This is impossible for $x \neq y$ since $v^1, v^2$ are linearly independent.

Secondly, if $s \in (0,1)$, then one can subtract $s(v^2)^\top x= s(v^2)^\top y$ from $((1-s)v^1 + sv^2)^\top x = ((1-s)v^1 + sv^2)^\top y$ and again get the two equations $(v^1)^\top x= (v^1)^\top y$, $(v^2)^\top x= (v^2)^\top y$. Thus, one ends up with the same contradiction as before.
\pend

\medskip
Next, find $k_\ell$ such that $\#X^<(w, p) + k_\ell = K$ and re-arrange the permutation \eqref{EqRotPermutation0} as follows
\begin{multline*}
x^{\pi_w(1)}, \ldots, x^{\pi_w(\#X^<(w^n, p))}, x^{\pi_{v^2}(k_1)}, \ldots, x^{\pi_{v^2}(k_\ell)}, \\ x^{\pi_{v^2}(k_\ell+1)}, \ldots, x^{\pi_{v^2}(k_L)}, x^{\pi_w(\#X^<(w^n, p)+L+1)}, \ldots, x^{\pi_w(N)}.
\end{multline*}
Re-label this permutation by $\pi_w$, its $K$-th element now is $x^{\pi_w(K)} = x^{\pi_{v^2}(k_\ell)}$. 

\medskip
{\bf Case (1) and (2).} In both cases, the next step is a rotation of $w$ in direction $v^2$: Set $w(s) = (1-s)w + sv^2$ and solve the problem
\[
\tag{RP} \text{maximize} \quad s
\]
subject to
\begin{eqnarray}
\label{EqRP1}
w(s)^\top x^{\pi_{w}(k)} & \leq & w(s)^\top x^{\pi_{w}(K)} \quad \text{for} \; k = 1, \ldots, K-1 \\
\label{EqRP2}
w(s)^\top x^{\pi_{w}(k)} & \geq & w(s)^\top x^{\pi_{w}(K)} \quad \text{for} \; k = K+1, \ldots, N \\
\label{EqRP3}
s & \in & [0, 1]  
\end{eqnarray}
The output is $w(\bar s)$ with some $\bar s \in (0, 1]$ and $X^<(w(\bar s), p)$, $X^=(w(\bar s), p)$. If $w(\bar s) \neq v^2$ (i.e., $\bar s \neq 1$), then $\#X^=(w(\bar s), p) \geq 2$. The idea is to keep $x^{\pi_w(K)}$ on the boundary of the $w(s)$-quantile for $s \in [0, \bar s]$.

\begin{lemma}[the rotation lemma for quantiles]
\label{LemRotationQ}
One has 

(a) $Q^-_{\tilde X, w(s)}(p ) = x^{\pi_w(K)} + H^+(w(s))$ for all $0 \leq s \leq \bar s$,

(b) strict inequalities in \eqref{EqRP1}, i.e., $X^<(w, p) \subseteq X^<(w(s), p)$ and $x^{\pi_w(K- k_\ell)}, \ldots, x^{\pi_w(K- 1)} \in X^<(w(s), p)$ for all $0< s < \bar s$,

(c)  strict inequalities in \eqref{EqRP2} for all $0< s < \bar s$,

(d)  $\#X^=(w(s), p) = 1$ for all $0< s < \bar s$.
\end{lemma}

{\sc Proof.} The first claim is by construction, i.e., \eqref{EqRP1}, \eqref{EqRP2}. If one assumes $w(t)^\top x^{\pi_{w}(n)} = w(t)^\top x^{\pi_{w}(K)}$ for some $t \in (0, \bar s)$ and $n \in \{1, \ldots, K-1\}$, then
\begin{equation}
\label{EqEqual1}
\of{(1-t)w + tv^2}^\top x^{\pi_{w}(n)} = \of{(1-t)w + tv^2}^\top x^{\pi_{w}(K)}.
\end{equation}
By \eqref{EqRotPermutation1}, one has
\[
w^\top x^{\pi_{w}(n)} \leq w^\top x^{\pi_{w}(K)}.
\]
If ``$=$" would be true in this inequality, then one had the two equations
\[
w^\top x^{\pi_{w}(n)} = w^\top x^{\pi_{w}(K)} \quad \text{and} \quad (v^2)^\top x^{\pi_{w}(n)} =  (v^2)^\top x^{\pi_{w}(K)}
\]
where the second is a consequence of the first and \eqref{EqEqual1}. This is a contradiction since the two vectors $w, v^2$ are linearly independent. 

If $w^\top x^{\pi_{w}(n)} < w^\top x^{\pi_{w}(K)}$ would be true, then \eqref{EqEqual1} would yield $(v^2)^\top x^{\pi_{w}(n)} > (v^2)^\top x^{\pi_{w}(K)}$. However, this cannot be true as one can see as follows. Taking $\eps > 0$ such that $t + \eps < \bar s$ one has by \eqref{EqRP1}
\[
\of{(1-t - \eps)w + (t + \eps)v^2}^\top x^{\pi_{w}(n)} = \of{(1-t-\eps)w + (t + \eps)v^2}^\top x^{\pi_{w}(K)}
\] 
\eqref{EqEqual1} would now yield
\[
(-\eps w + \eps v^2)^\top x^{\pi_{w}(n)} \leq (-\eps w + \eps v^2)^\top x^{\pi_{w}(K)}.
\]
Rearranging terms, using $w^\top x^{\pi_{w}(n)} \leq w^\top x^{\pi_{w}(K)}$ and dividing by $\eps > 0$ one would arrive at
\[
(v^2)^\top x^{\pi_{w}(n)} \leq (v^2)^\top x^{\pi_{w}(K)}
\]
which would produce a contradiction. \pend

\subsection{The quantile algorithm}
\label{SubSecQuantileAlg}

The proposed algorithm works as follows.

\medskip
{\bf Quantile Algorithm.}

\medskip
{\bf Step 1.} Initialize $w^1 := v^1$, $W := \{w^1\}$, $RI = \emptyset$.

\medskip
{\bf Step 2.} Until $w^n = v^2$ repeat: rotation step with input $w^n$, output $w^{n+1} := w(\bar s)$ and update $W := W \cup\{w^{n+1}\}$, update $n := n+1$. If $\#X^<(w^n, p) + \#X^=(w^n, p) = K$, then $RI := RI \cup\{n\}$.

\medskip
{\bf Step 3.} Update $W := W\bs\{w^k \mid k \in RI\}$. 

\medskip
{\bf Step 4.} Compute 
\begin{equation}
\label{EqAlgQuantile}
Q^-_{\tilde X, C}(p ) = \bigcap_{w \in W} \sqb{x(w, p) + H^+(w)}
\end{equation}

while removing the redundant $w$'s from $W$. 

\begin{corollary}
\label{CorTermination}
(1) The Quantile Algorithm terminates after a finite number of rotation steps. 

(2) \eqref{EqAlgQuantile} is true.
\end{corollary} 

{\sc Proof.} (1) There are only finitely many data points and hence only finitely many halfspaces in $\R^2$ with two points on their boundaries. The algorithm checks those which satisfy the conditions of Proposition \ref{PropIntNormals} consecutively which is guaranteed by Lemma \ref{LemRotationQ}, i.e., it identifies the set $W(p ) \cup\{v^1, v^2\}$.

(2) By construction, the algorithm identifies all elements of $W(p )$  and only removes those (in Step 2) which are redundant according to Proposition \ref{PropIntNormals}. \pend

\begin{remark}
If $\lceil Np\rceil = 1$, then $Q^-_{\tilde X, w}(p )$ is the convex hull of the data points plus the cone $C$.
\end{remark}

\begin{remark}
In our version of the algorithm, Step 4 makes use of the so-called Benson algorithm for the representation of convex polyhedrons implemented in the {\sc Bensolve} package \cite{LoehneWeissing17EJOR, BENSOLVE}. {\sc Bensolve} generates the smallest $H$-representation as well as the $V$-representation of $Q^-_{\tilde X, C}(p )$ which is a convex polyhedron. In particular, it removes the remaining redundant $w$'s. In the bivariate case, a more direct approach is possible, but we preferred to use {\sc Bensolve} since it is also usable for higher dimensional data.
\end{remark}

We show the result of the algorithm for the 20-data points example from above.

\begin{figure}[H]
\centering
\includegraphics[width=0.7\textwidth]{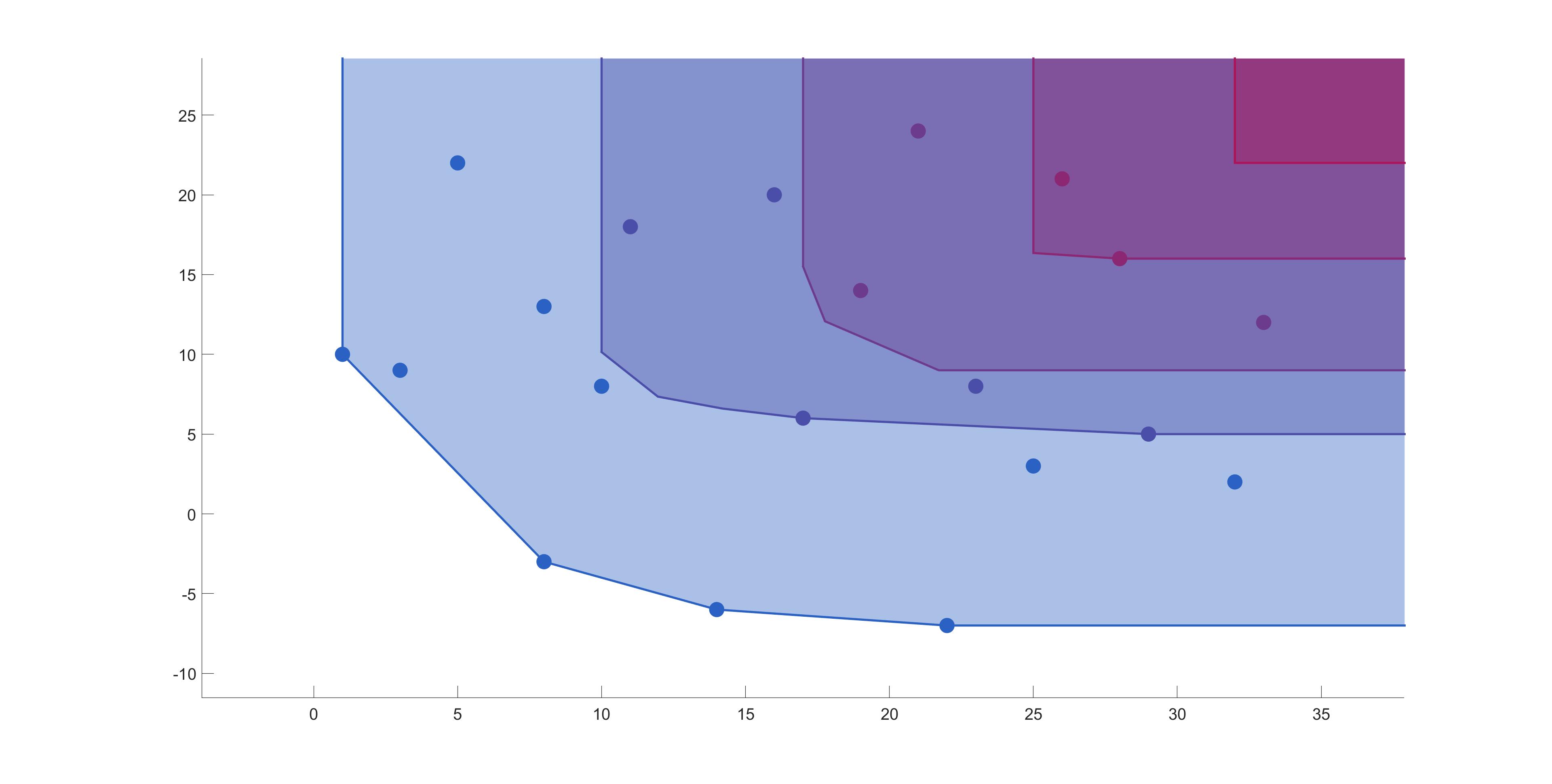}
\caption{$C = \R^2_+, N = 20, p = \{0.05,  0.275, 0.5, 0.725, 0.95\}$}
\label{fig:ConeQuantileLayer}
\end{figure}

\begin{example}
The figure \ref{fig:ConeQuantileLayer} illustrates the Quantiles for $p = \{0.05,  0.275, 0.5, 0.725, 0.95\}$ based on the 20-data points example from above.
\end{example}

In the worst case, each data point but two is on the boundary of two halfspaces which contribute to $Q^-_{\tilde X, C}(p )$.

\begin{example}
\label{ExWorstCase}
The following example showcases the situation in which indeed $N$ rotation steps have to be performed and the quantile set is the intersection of $N+1$ halfspaces. One may suspect that this feature is mainly due to the fact that the data points are not comparable with respect to $\leq_{\R^2_+}$.
\end{example}
\begin{figure}[H]
\centering
\includegraphics[width=0.7\textwidth]{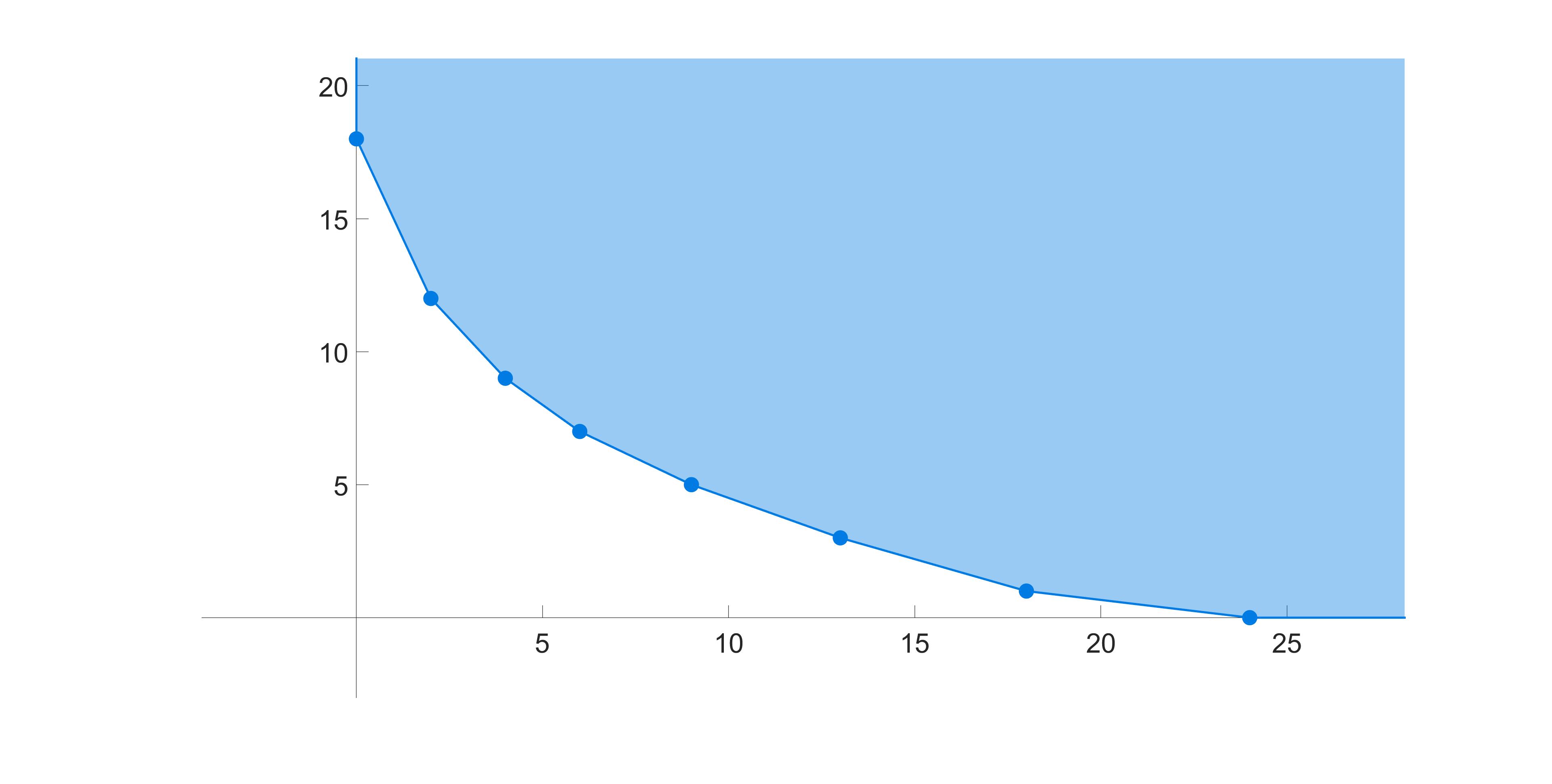}
\caption{$C = \R^2_+, N = 8, p = 0.125$}
\label{fig:ConeQuantileNsteps}
\end{figure}

\begin{remark}
\label{RemComparablePairs} If there are two data points $x, y \in \tilde X$ with $x - y \in \Int C$, then these two points cannot be on the boundary of a $w$-quantile for $w \in C^+\{0\}$. Indeed, if $v \perp (x-y)$, then $v \not\in C^+\bs\{0\}$ since $x - y \in \Int C$. This confirms that the quantile $Q^-_{\tilde X, C}(p )$ tends to have less vertices if there are more pairs of data points comparable with respect to $\leq_C$. In the extreme case of a linearly ordered data set, the lower $C$-quantile has the form "data point plus cone." Of course, the univariate case can be seen as a very special one for this situation.
\end{remark}

\subsection{The cone location depth algorithm}
\label{SubSecCLDepthAlg}

The algorithm in this subsection produces the value of $F_{\tilde X, C}(z)$ and $c\text{-}depth(z; \tilde X)$, respectively, for $z \in \R^d$. Basic ideas from the algorithm in Section \ref{SubSecQuantileAlg} will reappear, but the two algorithms are in some sense dual to each other: while the quantile algorithm changes the points on the boundary of the intermediate $w$-quantiles and keeps the property of being a $p$-quantile, the cone distribution function algorithm keeps the point $z \in \R^2$ at which $F_{\tilde X, C}$ is computed  on the boundary of intermediate $w$-quantiles, but the value of $w\text{-}depth(z; \tilde X)$ changes.

First, if necessary, $z \in \R^d$ is added to the set of data points $\tilde X$ and adjust $N$. This means, without loss of generality, $z \in \tilde X$ and $\# \tilde X = N$. For $w \in B^+$, let
\[
X^=(w, z) = \cb{x \in \tilde X \mid w^\top x = w^\top z}
\]
be the set of all data points which are on the boundary of the halfspace $z - H^+(w)$. Note the difference to $X^=(w, p)$ used previously.

Again, the algorithm is designed such that it starts with $w^1 = v^1$ and then runs through $B^+$ until it hits $v^2$. At an intermediate step, a $w^n \in B^+$ is generated. In each of the three cases 

(1) $w^n = v^1$ and $\#X^=(v^1, z) = 1$, i.e., $X^=(v^1, z) = \{z\}$,

(2) If $w^n \neq v^2$ and $\#X^=(w^n, z) \geq 2$,

(3) $w^n = v^2$,
\\
a value $K_n$ is computed which is a current upper bound of $c\text{-}depth(z; \tilde X)$ and then a rotation step is carried out except in case (3) which again serves as a stopping criterion. In case (1) and (2), the rotation step is performed such that $z$ stays on the boundary of the resulting $w^{n+1}$-quantile.

\medskip
{\bf Case (1).} Find the permutation
\begin{equation}
\label{EqRotPermutation00}
w^\top x^{\pi_w(1)} \leq \ldots < w^\top z < \ldots \leq w^\top x^{\pi_w(N)}
\end{equation}
and the number $K$ such that $x^{\pi_w(K)} = z$.

\medskip
{\bf Case (2).} The input is $w \in B^+\bs\{v^2\}$. Find the permutation
\begin{equation}
\label{EqRotPermutation2}
w^\top x^{\pi_w(1)} \leq \ldots \leq w^\top z \leq \ldots \leq w^\top x^{\pi_w(N)}
\end{equation}
Find the set $X^=(w, z)$ and determine $L := \#X^=(w, z)$ as well as $k := \#\cb{x \in \tilde X \mid w^\top x < w^\top z}$. Clearly, one always has $z \in X^=(w, z)$. On $X^=(w, z)$, define a new permutation $\pi_2$ by
\begin{equation}
\label{EqRotPermutationEqual}
(v^2)^\top x^{\pi_2(1)} \leq \ldots \leq (v^2)^\top z \leq \ldots \leq (v^2)^\top x^{\pi_2(L)}.
\end{equation}
As in Proposition \ref{PropStrictPermut}, all inequalities in this permutation are strict. Let $\ell$ be the number such that $z = x^{\pi_2(\ell)}$. Set $K = k + \ell$. In the permutation \eqref{EqRotPermutation2}, replace the elements in $X^=(w, z)$ in the order generated by the permutation \eqref{EqRotPermutationEqual} and relabel the permutation by $\pi_w$. Now, $z = x^{\pi_w(K)}$.

\medskip
{\bf Case (1) and (2).} In both cases, set $w(s) = (1-s)w + sv^2$ and solve the problem (RP), \eqref{EqRP1}-\eqref{EqRP3}. The output is $w(\bar s)$ with some $\bar s \in (0, 1]$. If $w(\bar s) \neq v^2$ (i.e., $\bar s \neq 1$), then $\#X^=(w(\bar s), z) \geq 2$.

\begin{lemma}[the rotation lemma for the CDF]
\label{LemRotationF} If $w \in B^+\bs\{v^2\}$, then

(a) there are $k+L$ data points in $z - H^+(w)$ (with $L =1$ in Case (1)), i.e., $w\text{-}depth(z; \tilde X) = k + L$,

(b) there are $k+\ell$ data points in $z - H^+(w(s))$ for all $s \in (0, \bar s)$, i.e., $w(s)\text{-}depth(z; \tilde X) = k + \ell$ for $s \in (0, \bar s)$,

(c) there are at least $k+\ell$ data points in $z - H^+(w(\bar s))$, i.e., $w(\bar s)\text{-}depth(z; \tilde X) \geq k + \ell$.
\end{lemma}

{\sc Proof.} First, observe that one has $X^=(w(s), z) = \{z\}$ for $s \in (0, \bar s)$ due to the fact that there are only finitely many data points: the rotation of $w$ in the direction of $v^2$ around $z$ removes all data points from the boundary line of $z - H^+(w)$ except $z$.

(a) This is due to the definition of $k$ and $L$. (b) According to \eqref{EqRotPermutationEqual} and \eqref{EqRP2}, the $L-\ell$ data points $x^{\pi_2(\ell+1)}, \ldots, x^{\pi_2(L)}$ are not in $z - H^+(w(s))$ along with the data points $x^{\pi_2(L+1)}, \ldots, x^{\pi_2(N)}$, while the points $x^{\pi_2(1)}, \ldots, x^{\pi_2(\ell)}$ remain in $z - H^+(w(s))$ according to \eqref{EqRP1} for all $s \in (0, \bar s)$. (c) According to \eqref{EqRP1}, \eqref{EqRP2}, there are at least $k + \ell$ data points in $z - H^+(w(\bar s))$ and if at least one of the inequalities in \eqref{EqRP2} is satisfied as an equation for $s = \bar s$, then there are at least $k+ \ell + 1$ data points in $z - H^+(w(\bar s))$. \pend

\medskip
{\bf The cone location depth algorithm.}

\medskip
{\bf Step 1.} Initialize $w^1 := v^1$ and set $K_1 := K$ with $K$ as found in Case (1) or (2).

\medskip
{\bf Step 2.} Until $w^n = v^2$ repeat: starting with $\pi_w$ obtained from the permutations found in Case (1) or (2), perform a rotation step with input $w^n$, output $w^{n+1} := w(\bar s_n)$ and $K_{n+1}$. Update $K :=\min\{K_n, K_{n+1}\}$. Update $n:=n+1$.

\medskip
{\bf Step 3.} If $w^n = v^2$, compute $(v^2)\text{-}depth(z; \tilde X)$ and update $K := \min\{K, (v^2)\text{-}depth(z; \tilde X)\}$.

\medskip
{\bf Step 4.} Compute 
\begin{equation}
\label{EqAlgCDF}
F_{\tilde X, C}(z) = 
\left\{
\begin{array}{ccc}
\frac{K}{N} & : & z \; \text{is an original data point} \\[.2cm]
\frac{K-1}{N-1} & : & z \; \text{is not an original data point} 
\end{array}
\right.
\end{equation}

\begin{corollary}
(1) The Cone Distribution Function algorithm terminates after a finite number of rotation steps.

(2) \eqref{EqAlgCDF} is true.
\end{corollary}

{\sc Proof.} (1) There are only finitely many data points and hence only finitely many halfspaces in $\R^2$ with at least two data points on their boundaries. The algorithms checks those with $z$ as one boundary point and normal direction in $B^+$ as well as $z - H^+(v^1)$, $z - H^+(v^2)$.

(2) By Lemma \ref{LemRotationF}, $K = k + \ell \leq w(s)\text{-}depth(z; \tilde X)$ for each $s \in [0, \bar s]$ with equality for $s \in (0, \bar s)$. By construction in Step 2, the algorithm determines the minimum of these numbers and $v^2\text{-}depth(z; \tilde X)$. Equation \eqref{EqAlgCDF} follows taking into account that $z$ might or might not be an original data point. \pend

\subsection{Pseudocodes}

In this section, pseudocodes for the two algorithms are provided along with a few explanatory remarks. The following algorithms assigns an index to each element of the set of data points $\tilde X$. This index is then used trough out the algorithm in order to simplify the identification of each point and to minimize rounding errors.

\begin{algorithm}[H]
\caption{Bivariate Lower Cone Quantile.}\label{algo_CQ1}
\begin{algorithmic}[1]
\Procedure{ConeQuantile}{$p$,$\tilde X$,$b$}\Comment{$p \in (0,1]$ , $\tilde X \in \R^N \times\R^2$, $b \in \R^2\times\R^2$}
\State $v \leftarrow \begin{pmatrix} b_{1,2} &  -b_{1,1} \\ -b_{2,2} &  b_{2,1} \end{pmatrix}$ \Comment{see remark \ref{remark_algo_v}}
\If {$(v^1)^Tb^2 < 0 \vee (v^2)^Tb^1 < 0$}
	\State $v \leftarrow \begin{pmatrix} -b_{1,2} &  b_{1,1} \\ b_{2,2} &  -b_{2,1} \end{pmatrix}$
\EndIf
\State $K \leftarrow \lceil p \times N \rceil$ \Comment{see proposition \ref{PropDataPtBoundary}}
\State $w^1 \leftarrow v^1$, $W \leftarrow \cb{w^1}$ \Comment{$v^1$ is assigned to $w^1$, which is then saved in $W$}
\State $n=1$ \Comment{$n$ counts the iterations}
\While {}
	\State $w \leftarrow w^n$
	\State $\pi^{\tilde X}_w \leftarrow$ \Call{IndexSort}{$\tilde X^\top w$} \Comment{see remark \ref{remark_algo_indexsort}}
	\State $q \leftarrow q \cup  \cb{{w^\top x^{\pi^{\tilde X}_w(K)}}}$ \Comment{see remark \ref{remark_algo_scalarization}}
	\State $B_R \leftarrow \cb{j \in \{1,\hdots , N\} \mid w^\top {x^{\pi^{\tilde X}_w(K)}} = w^\top x^{\pi^{\tilde X}_w(j)}}$ \Comment{see remark \ref{remark_algo_CD_4}}
	\State $B_I \leftarrow \cb{\pi^{\tilde X}_w \in \{1,\hdots , N\} \mid w^\top {x^{\pi^{\tilde X}_w(K)}} = w^\top x^{\pi^{\tilde X}_w}}$ \Comment{index of boundary points}
	\State $B \leftarrow \cb{x \in \tilde X \mid w^\top {x^{\pi^{\tilde X}_w(K)}} = w^\top x}$ \Comment{set of boundary points}
	\If {$\max(B_R)>K$} \Comment{see remark \ref{remark_algo_CQ_redundant}}
		\State $I_K \leftarrow I_K \cup n$
	\EndIf 	
	\If {$w=v^2$} \Comment{stopping condition, see remark \ref{remark_algo_CQ_stop}}
		\State {\bf break while}
	\EndIf
	\State $\pi_{v^2}^B \leftarrow$ \Call{IndexSort}{$B^\top v^2$}
	\State $\pi^{\tilde X}_w(B_R) \leftarrow B_I(\pi_{v^2}^B)$ \Comment{see remark \ref{remark_algo_rearrange}}
	\State $w^{n+1} \leftarrow$ \Call{Rotation}{$\tilde X, w, v^2, \pi^{\tilde X}_w, K$} \Comment{see (RP), \eqref{EqRP1}-\eqref{EqRP3}}
	\State $W \leftarrow W \cup \cb{w^{n+1}}$ \Comment{add the new direction $w^{n+1}$ to $W$}
	\State $n = n + 1$ \Comment{update the counting variable $n$}
\EndWhile
\algstore{algo_CQ}
\end{algorithmic}
\end{algorithm}

\begin{algorithm}[H]
\caption*{Bivariate Lower Cone Quantile.}\label{algo_CQ2}
\begin{algorithmic}[1]
\algrestore{algo_CQ}
\State $I_K \leftarrow I_K \cup \cb{1, n}$ \Comment{see remark \ref{remark_algo_CQ_redundant_v}}
\State $W \leftarrow W(I_K)$ \Comment{extract the directions with indexes $I_K$}
\State $q \leftarrow q(I_K)$ \Comment{extract the scalarized boundary points with indexes $I_K$}
\State $Q \leftarrow$ \Call{polyh}{$W,q$} \Comment{see remark \ref{remark_algo_CQ_BENSOLVE}}
\State $W,q \leftarrow$ \Call{hrep}{$Q$}
\State $D,V \leftarrow$ \Call{vrep}{$Q$}
\EndProcedure
\end{algorithmic}
\end{algorithm}

\begin{algorithm}[H]
\caption{Bivariate Cone Distribution function.}\label{algo_CD1}
\begin{algorithmic}[1]
\Procedure{ConeDistribution}{$z$,$\tilde X$,$b$}\Comment{$z \in \R^2$ , $\tilde X \in \R^N \times\R^2$, $b \in \R^2\times\R^2$}
\State $I \leftarrow \cb{1,\hdots,N}$ \Comment{index set of $\tilde X = \cb{x^1,\hdots,x^N}$}
\State $I_z \leftarrow \cb{i \in I \mid z = x^i}$ \Comment{set of positions of $z$ in $\tilde X$}
\If {$\#(I_z) = 0$} \Comment{see remark \ref{remark_algo_CD_1}}
	\State $\tilde X \leftarrow \tilde X \cup z$ \Comment{add $z$ at the end of $\tilde X$}
	\State $N \leftarrow N+1$	
	\State $i_z \leftarrow N$
	\State $I \leftarrow I \cup i_z$
\ElsIf {$\#(I_z) = 1$}
	\State $i_z \leftarrow I_z$ \Comment{index of $z$ in $\tilde X$}
\ElsIf {$\#(I_z) > 1$}
	\State $i_z \leftarrow max(I_z)$ \Comment{see remark \ref{remark_algo_CD_2}}
\EndIf
\State $v \leftarrow \begin{pmatrix} b_{1,2} &  -b_{1,1} \\ -b_{2,2} &  b_{2,1} \end{pmatrix}$ \Comment{see remark \ref{remark_algo_v}} 
\If {$(v^1)^Tb^2 < 0 \vee (v^2)^Tb^1 < 0$}
	\State $v \leftarrow \begin{pmatrix} -b_{1,2} &  b_{1,1} \\ b_{2,2} &  -b_{2,1} \end{pmatrix}$
\EndIf
\algstore{algo_CD}
\end{algorithmic}
\end{algorithm}

\begin{algorithm}[H]
\caption*{Bivariate Cone Distribution function.}\label{algo_CD2}
\begin{algorithmic}[1]
\algrestore{algo_CD}
\State $w \leftarrow v^1$, $k \leftarrow N$ \Comment{$w$ and $k$ are updated in each iteration}
\While {}
	\State $\pi^{\tilde X}_w \leftarrow$ \Call{IndexSort}{$\tilde X^\top w$} \Comment{see remark \ref{remark_algo_indexsort}}
	\State $B_R \leftarrow \cb{j \in \{1,\hdots , N\} \mid w^\top x^{i_z} = w^\top x^{\pi^{\tilde X}_w(j)}}$ \Comment{see remark \ref{remark_algo_CD_4}}
	\If {$w=v^2$} \Comment{stopping condition, see remark \ref{remark_algo_CD_stop}}
		\State $k \leftarrow min(k,max(B_R))$
		\State {\bf break while}
	\EndIf
	\State $B_I \leftarrow \cb{\pi^{\tilde X}_w \in \{1,\hdots , N\} \mid w^\top x^{i_z} = w^\top x^{\pi^{\tilde X}_w}}$ \Comment{Index of boundary points}
	\State $B \leftarrow \cb{x \in \tilde X \mid w^\top x^{i_z} = w^\top x}$ \Comment{set of boundary points}
	\State $\pi_{v^2}^B \leftarrow$ \Call{IndexSort}{$B^\top v^2$}
	\State $\pi^{\tilde X}_w(B_R) \leftarrow B_I(\pi_{v^2}^B)$ \Comment{see remark \ref{remark_algo_rearrange}}
	\State $K \leftarrow \cb{j \in \{1,\hdots , N\} \mid i_z = \pi^{\tilde X}_w(j)}$ \Comment{see remark \ref{remark_algo_CD_position_iz}} 
	\State $k \leftarrow min(k,K)$ \Comment{update $k$ if $K < k$}
	\State $w \leftarrow$ \Call{Rotation}{$\tilde X, w, v^2, \pi^{\tilde X}_w, K$} \Comment{see (RP), \eqref{EqRP1}-\eqref{EqRP3}}
\EndWhile
\If {{$\#(I_z) = 0$}} \Comment{$p$ is equal to $\frac{k}{N}$ or $\frac{k-1}{N-1}$, if $z \notin \tilde X$.}
	\State $p \leftarrow \frac{k-1}{N-1}$ 
\Else
	\State $p \leftarrow \frac{k}{N}$
\EndIf
\EndProcedure
\end{algorithmic}
\end{algorithm}

\begin{remark}
\label{remark_algo_CD_1}
$I_z$ is the set of positions of $z$ in $\tilde X$, whereas $\#$ is a counting function. There are three important cases: $z \notin \tilde X$, $\#(I_z) = 1$ and $\#(I_z) > 1$.
\end{remark}

\begin{remark}
\label{remark_algo_CD_2}
As the sort algorithm will keep the order of the multiples of $z$, we need the index ($i \in I$) of the ``last'' multiple in $\tilde X$.
\end{remark}

\begin{remark}
\label{remark_algo_v}
The V-representation of $C$ (matrix $b$) is converted into the V-representation of $C^+$ (matrix $v$). $v^1$ and $v^2$ are chosen orthogonal to $b^1$ and $b^2$, respectively, such that $(v^1)^Tb^2 \geq 0$, $(v^2)^Tb^1 \geq 0$.
\end{remark}

\begin{remark}
\label{remark_algo_indexsort}
The function $indexsort(B)$ sorts the elements of $B$ in ascending order and outputs the indexes of the elements of B in sorted order.
\end{remark}

\begin{remark}
\label{remark_algo_CD_4}
The rank/position $j$ of the boundary points in the permutation $\pi_w^{\tilde X}(j)$.
\end{remark}

\begin{remark}
\label{remark_algo_CD_stop}
The while loop is stopped, if the $w$ generated in the previous iteration is equal to $v_2$. Moreover, $k$ is updated if $(\#X^<(v^2, p) + \#X^=(v^2, p)) < k$.
\end{remark}

\begin{remark}
\label{remark_algo_rearrange}
Rearrange elements of $B_I$ based on the permutation $\pi^{B}_{v^2}$ and replace this rearranged indexes with the ones in the permutation $\pi^{\tilde X}_{w}$ at the positions $B_R$.
\end{remark}

\begin{remark}
\label{remark_algo_CD_position_iz}
The $K$ for each iteration (direction) is found by looking for the position of the index of $z$ ($i_z$) in the permutation $\pi^{\tilde X}_w$.
\end{remark}

\begin{remark}
\label{remark_algo_scalarization}
The point $x^i \in \tilde X$ with index $i = \pi^{\tilde X}_w(K)$ is scalarized with $w$ and saved to $q$.
\end{remark}

\begin{remark}
\label{remark_algo_CQ_redundant}
The halfspaces that include more than K points are relevant for the calculation of the quantile. The number of points in the halfspace can be derived by taken the biggest number in the set $B_R$, which includes the ranks $j$ of the boundary points in the permutation $\pi^{\tilde X}_w(j)$. 
\end{remark}

\begin{remark}
\label{remark_algo_CQ_stop}
The while loop is stopped, if the $w$ generated in the previous iteration is equal to $v_2$.
\end{remark}

\begin{remark}
\label{remark_algo_CQ_redundant_v}
The halfspaces that include K points are redundant, with exception of $v^1$ and $v^2$. Therefore, the index $1$ and $n$ are also included to $I_K$.
\end{remark}

\begin{remark}
\label{remark_algo_CQ_BENSOLVE}
The {\sc Bensolve} tools (see \cite{{LoehneWeissing17EJOR}, BENSOLVE}) provide two functions, {\sc polyh} and {\sc vrep}, to convert a H-representation into a V-representation. First, the H-representation is transformed into the 'polyhedral object' format with {\sc polyh}. Secondly, the output of the latter function is then used as input for {\sc vrep} which returns the V-representation. {\sc Bensolve} is also used to remove all redundant halfspaces via {\sc hrep}.
\end{remark}

\subsection{Computing Tukey depth regions and Tukey depth functions}

The same algorithm can be used to compute bivariate Tukey depth regions. In this case, it is applied consecutively to the three sets
\begin{align*}
B^+_0 & = \co\cb{(-1, -1)^\top, (1, -1)^\top}, \\ 
B^+_1 & = \co\cb{(1, -1)^\top, (0, 1)^\top}, \\ 
B^+_2 & = \co\cb{(-1, -1)^\top, (0, 1)^\top}.
\end{align*}
The resulting algorithm is very close to the one described in \cite{RutsRousseeuw96CSDA}. It is worth noting though that replacing the unit circle used in \cite{RutsRousseeuw96CSDA} by the three sections $B^+_0, B^+_1, B^+_2$ admits the formulation of the rotation step as a linear problem. The following picture shows a Tukey depth region for Example \ref{ExDepth} with 20 data points.

\begin{figure}[H]
  \centering
  \includegraphics[width=1\textwidth]{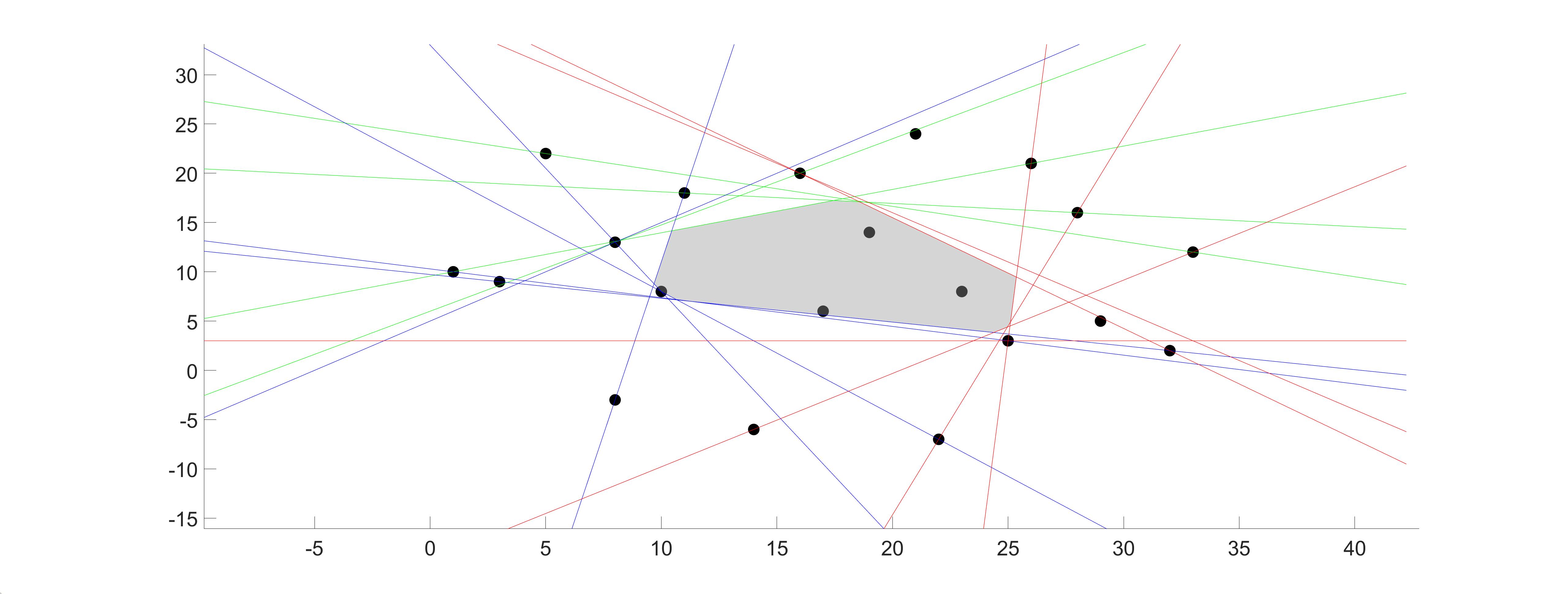}
  \caption*{Tukey depth region for $p=0.25$. Halfspace boundaries: {\color{green}$B^+_0$}, {\color{blue}$B^+_1$}, {\color{red}$B^+_2$}}
  \label{fig:ConeQuantileTukey}
\end{figure}

The cone location depth algorithm was applied to compute the values of the halfspace depth function in Example \ref{ExDepth} above.

\section{Examples}
\label{SecEx}

\begin{example}
\label{ExNormal}
The following figure shows the values of the cone location depth for a sample of the bivariate standard normal distribution which comprises 80 points with $C=\R^2_+$.
\end{example}

\begin{figure}[H]
	\centering
	\includegraphics[width=.9\textwidth]{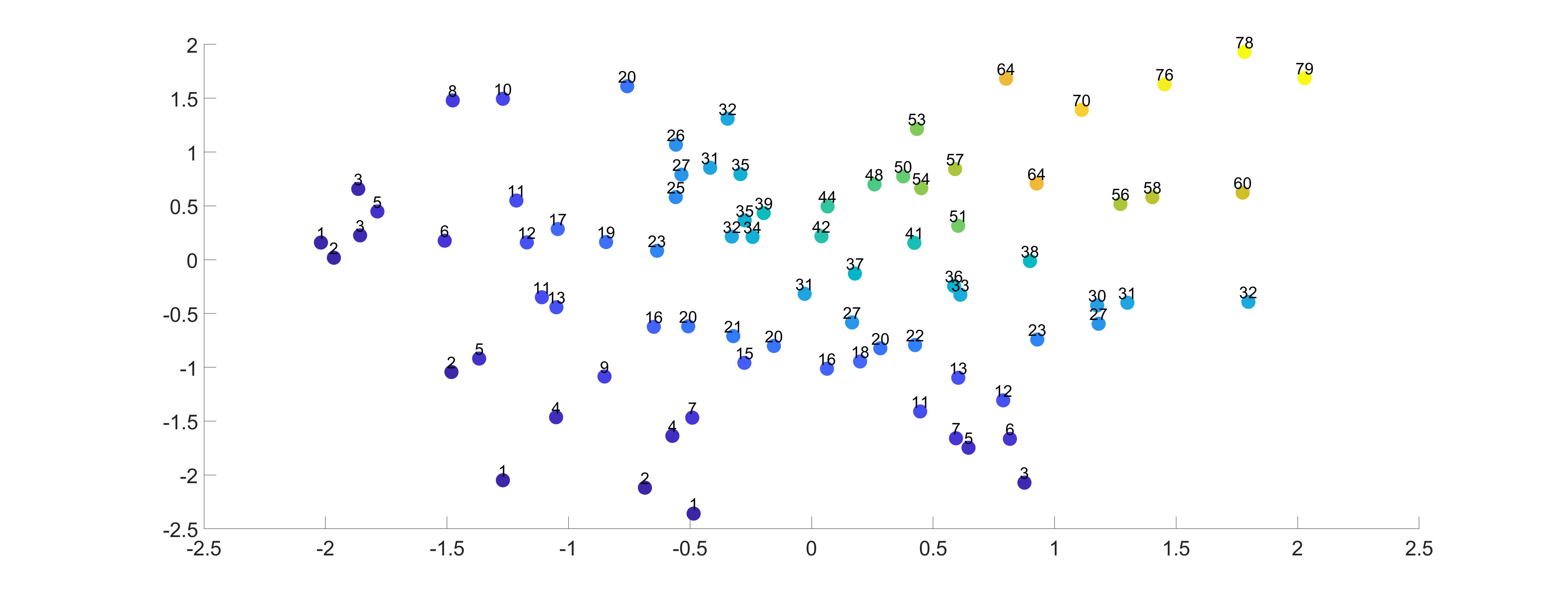}
	\caption{Cone location depth for a sample of a bivariate normal distribution}
	\label{fig:TwoNorm}
\end{figure}

\begin{example}
\label{ExUniformSquares} The next figure shows the values of the cone location depth for a sample of a bivariate uniform distribution over two squares with the cone $C=\R^2_+$. This example is taken from \cite{BelloniWinkler11AS} and was also discussed in \cite{HamelKostner18JMVA}.
\end{example}
\begin{figure}[H]
	\centering
	\includegraphics[width=.8\textwidth]{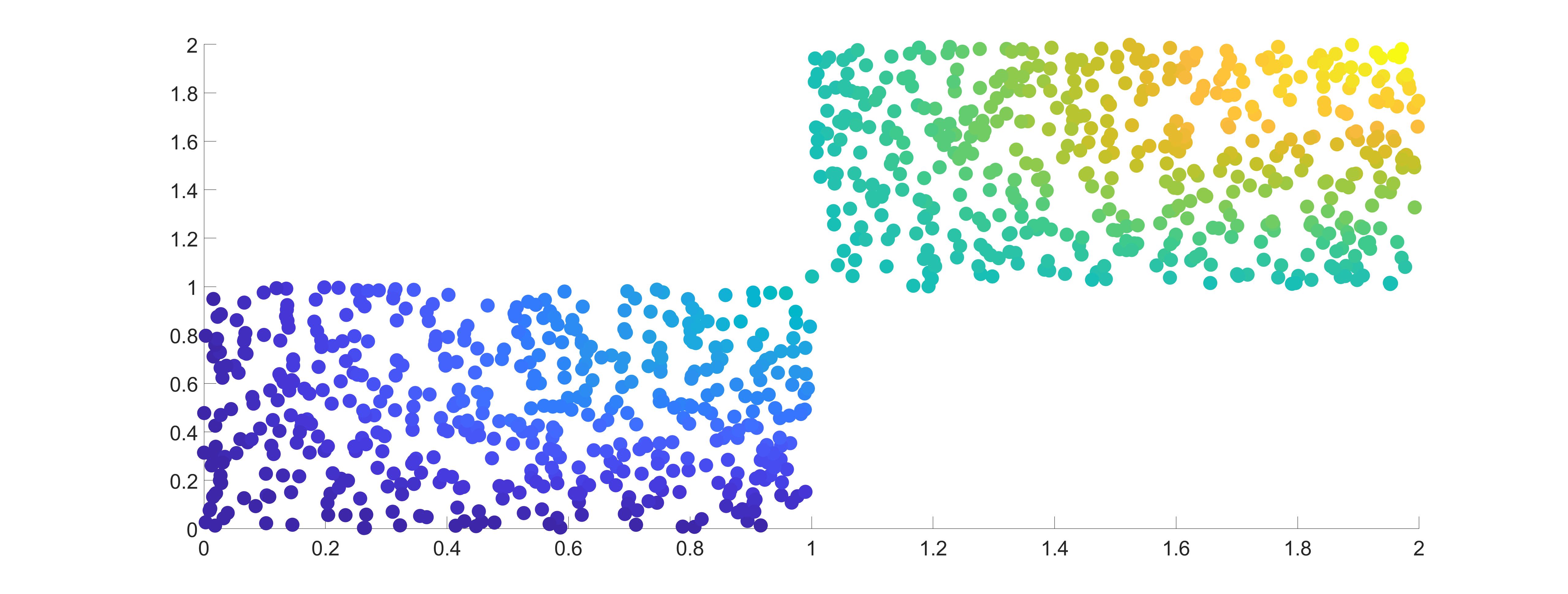}
	\caption{Cone location depth for a sample of $X \sim Uniform\big[(0,1)\times(0,1)\cup(1,2)\times(1,2)\big]$.}
	\label{fig:TwoNorm}
\end{figure}

A few examples with real world data are added. The aim is not to provide an extensive analysis for them, but to show how a bivariate approach might change the picture compared to a univariate one. One should note that in each of these examples there is an intuitive understanding for "better" or "worse" data points. Consequently, a mere depth function approach would not produce meaningful results. For all of these examples, the cone $\R^2_+$ generating the componentwise order is used for the sake of simplicity, but there are of course other (and maybe better) options.

\begin{example}
\label{ExHurricanes} The motivation to consider this data set stems from \cite{HebertWeinzapfelChambers10, KlotzbachEtAl20BAMS} where it is discussed that the maximum sustained wind speed (WS) alone (and hence the storm category according to the Saffir-Simpson scale) is not always a good proxy for the potential destructive impact of a hurricane. While the minimal sea level pressure (SLP) already seems to be a better proxy (see the strong arguments in \cite{KlotzbachEtAl20BAMS}), one can also think of using more than one parameter: along with wind speed and central pressure, storm wind radii and the storm translation speed (in particular, the devastating impact of Hurricane Dorian 2019 on the Bahamas supports such a parameter) are discussed in the quoted references. Yet another parameter should be storm surge as discussed in \cite{WalkerEtAl18WF} since this feature was removed from the Saffir-Simpson scale in the wake of 2005 Katrina.

Here, the two parameters WS and SLP are used to compare hurricanes. Note that the values usually occur at different times which are also different from the time of landfall: other choices are possible and can easily be implemented. The cone $C$ is generated by $\cb{b^1, b^2} =\cb{(1, 0), (0, -1)}$, since stronger wind and lower pressure characterize potentially more destructive hurricanes. The quantiles could then be used to categorize hurricanes.

A comparison with a scale based only on SLP produced very similar results, but with the following difference: The $C$-quantiles categorize the hurricanes based on the SLP as well as WS. This means that a hurricane needs to overcome a threshold in both variables. Therefore, a hurricane, that  is in a specific category on a scale based on one variable is not always in the corresponding category on the scale based on $C$-quantiles, as it does not reach the threshold in the other variable.

This serves just as an example to illustrate the potential of a categorization via set-valued quantiles. We think that even different cones should be used. The one generated by $b^1 = (1, 0)$ and $b^2 = \sqrt{2}/2(-1, -1)$ (or similar vectors) would be a good candidate to capture destructive storms with very low central pressure, but only moderate wind speed: the 2012 hurricane Sandy is the prominent example since it was not even categorized as a hurricane at landfall according to the Saffir-Simpson scale. On the other hand, the half-space $H^+((0, -1))$ as cone $C$ reproduces the situation considered in \cite{KlotzbachEtAl20BAMS}. The question how to choose the cone should be subject to an extended analysis using historical data. 

Such data can also be utilized in the following way: the track of a current hurricane can be followed in the quantile graph, thus providing a strong impression how "close" it comes to previous (major) hurricane (see Dorian's 2019 track below). This could deliver a strong warning message.

Finally, note that fast computations in real-time (see \cite{KlotzbachEtAl20BAMS}) are an issue. The algorithms presented in this paper could meet such demands.
\begin{figure}[H]
	\centering
	\includegraphics[width=1\textwidth]{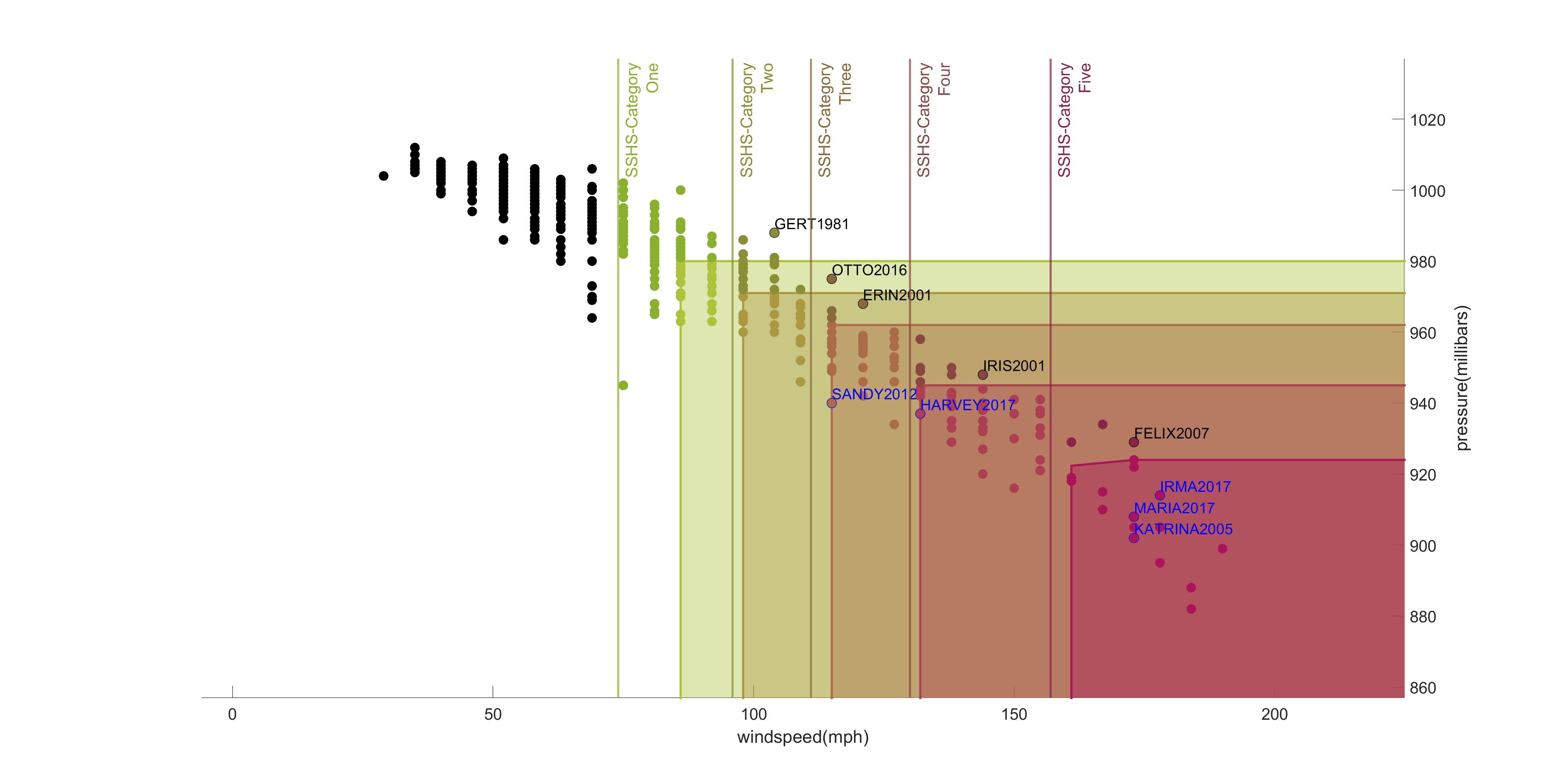}
	\caption{The figure shows 5 categories based on the $C$-quantiles.}
	\label{fig:HurricaneMinPres_categories}
\end{figure}

\begin{figure}[H]
	\centering
	\includegraphics[width=1\textwidth]{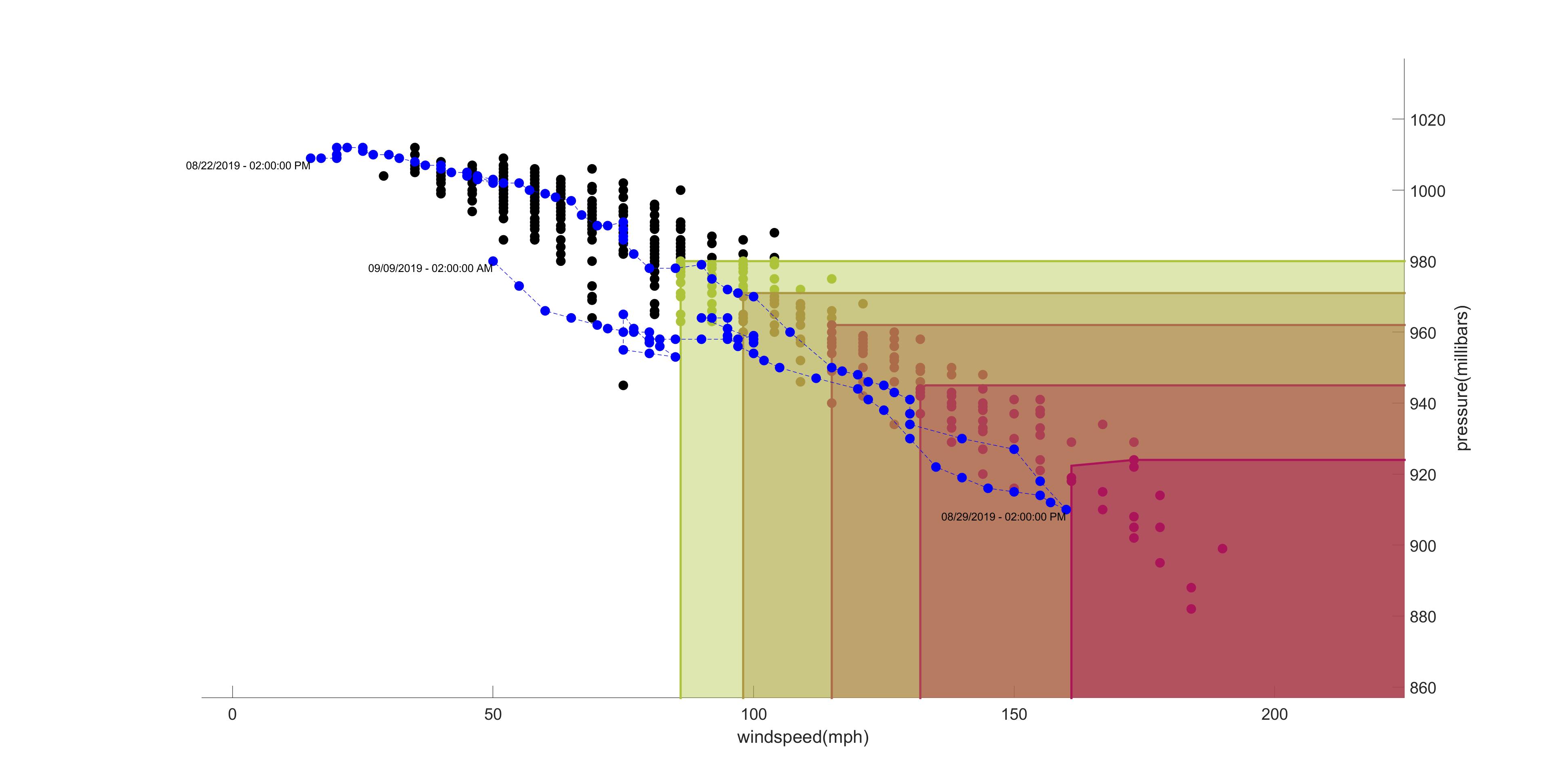}
	\caption{Dorian's 2019 track through the quantile graph.}
	\label{fig:Dorian}
\end{figure}
\end{example}

\begin{example}
South Tyrol, a region in the north-east of Italy, is one of the biggest apple producer in Europe. Due to its geographical position in the middle of the Alps it is hit regularly by severe hailstorms. The three main possibilities to hedge against the risk of hail damages are a hail insurance with substantial European subsidies, the installation of hail nets and the combination of these two options. In the following analysis the third option is compared to the first one. Hence, two types of farmers are compared: one type has stipulated only a hail insurance and the other one hedges their product with hail nets and insures it additionally. Both insurance products are subsidized by the European Union, whereas the contract with hail nets is discounted.

The data used for this analysis is given by the Hagelschutzkonsortium. It comprises, inter alia, the area insured ($AI$), the sum insured ($SU$), the premia paid by the farmer ($FP$), the sustained damage ($SD$) and the indemnity payments ($IP$) for each insurance contract signed between 2013 to 2017 by the members of the hail-defense syndicate. Only the farmers that have signed a contract in all 5 years are taken for the analysis. Each farmer is represented by two numbers. The first number is a proxy for the yearly business return: 
\begin{align*}
\text{no hail net: } & RB=SU-SD \\
\text{with hail net: } & RB=SU-SD-1652.98*AI 
\end{align*}
where the estimated cost for a hectar of hail net is \EUR{1652.98} per year, for detailed information see \cite{BallatoreVittone2008}.
The second number is an indicator if the insurance pays off:
\[
RI=IP-FP
\]

The figure \ref{fig:hail_data} has on the x-axes $RB$-values and on the y-axes $RI$-values. The blue data points represent the farmers with no hail nets and that have signed an insurance contract in each year between 2013 and 2017, whereas the red data points represent the farmers with hail nets. There are points that are located far away from the rest, this actually corresponds to the market situation with a few very big producers and the vast majority of small farmers.

\begin{figure}[H]
	\centering
	\includegraphics[width=0.8\textwidth]{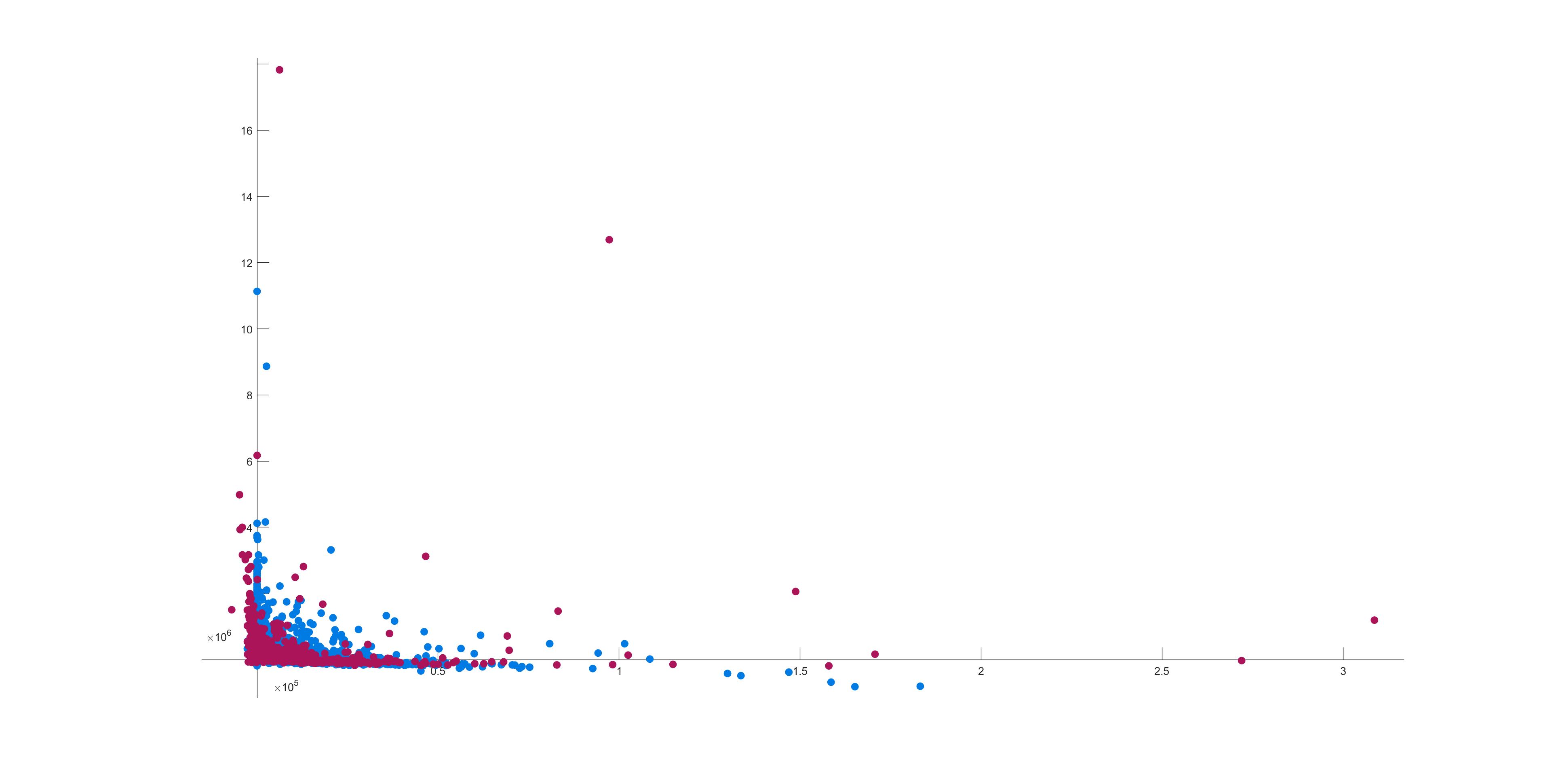}
	\caption{The farmers with hail nets are the red points.}
	\label{fig:hail_data}
\end{figure}

\begin{figure}[H]
\label{fig:hailQ}
\centering
\begin{minipage}{.31\textwidth}
  \centering
  \includegraphics[width=1\textwidth]{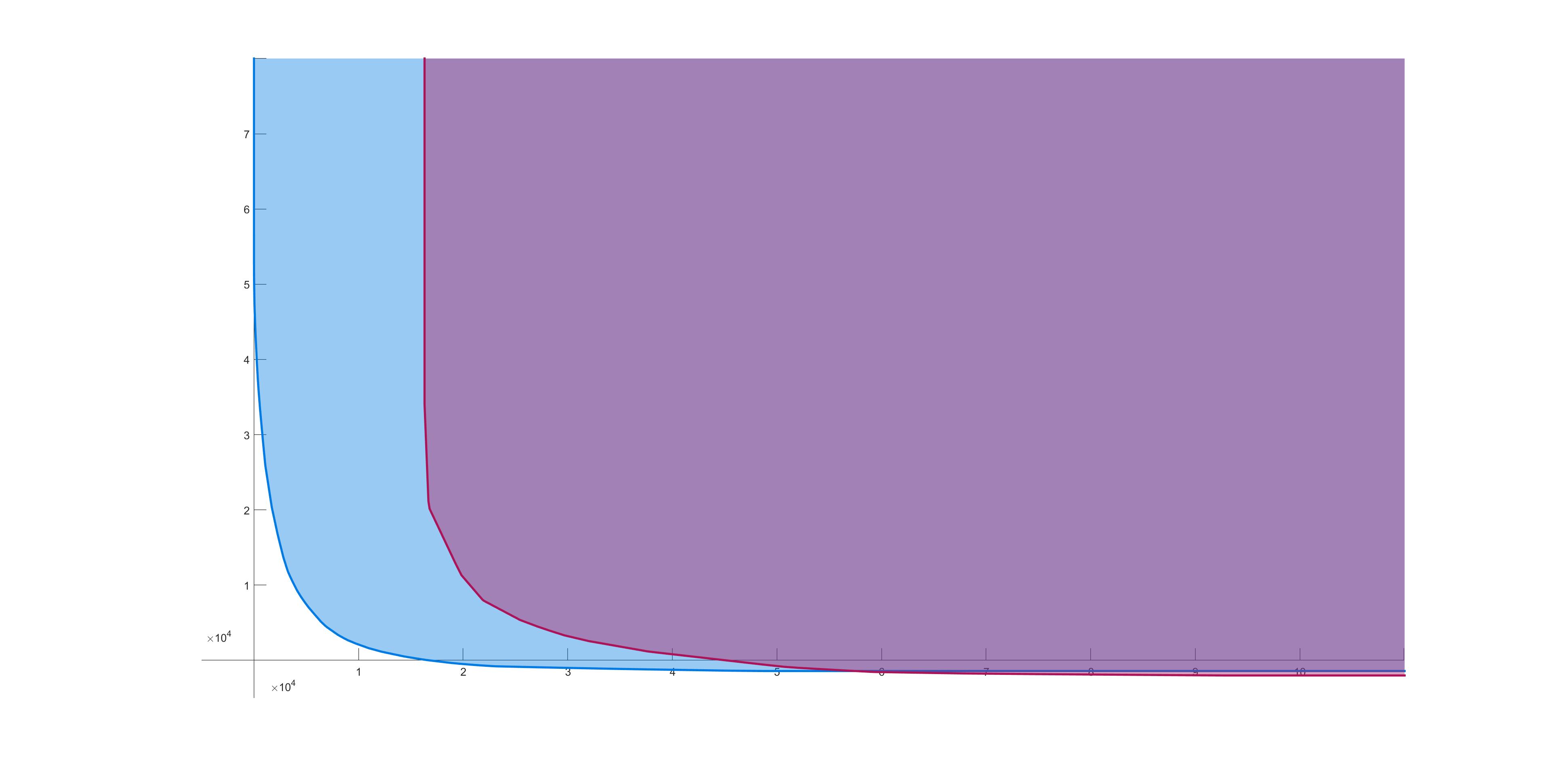}
  \caption*{$p=0.25$}
\end{minipage}
\begin{minipage}{.31\textwidth}
  \centering
  \includegraphics[width=1\textwidth]{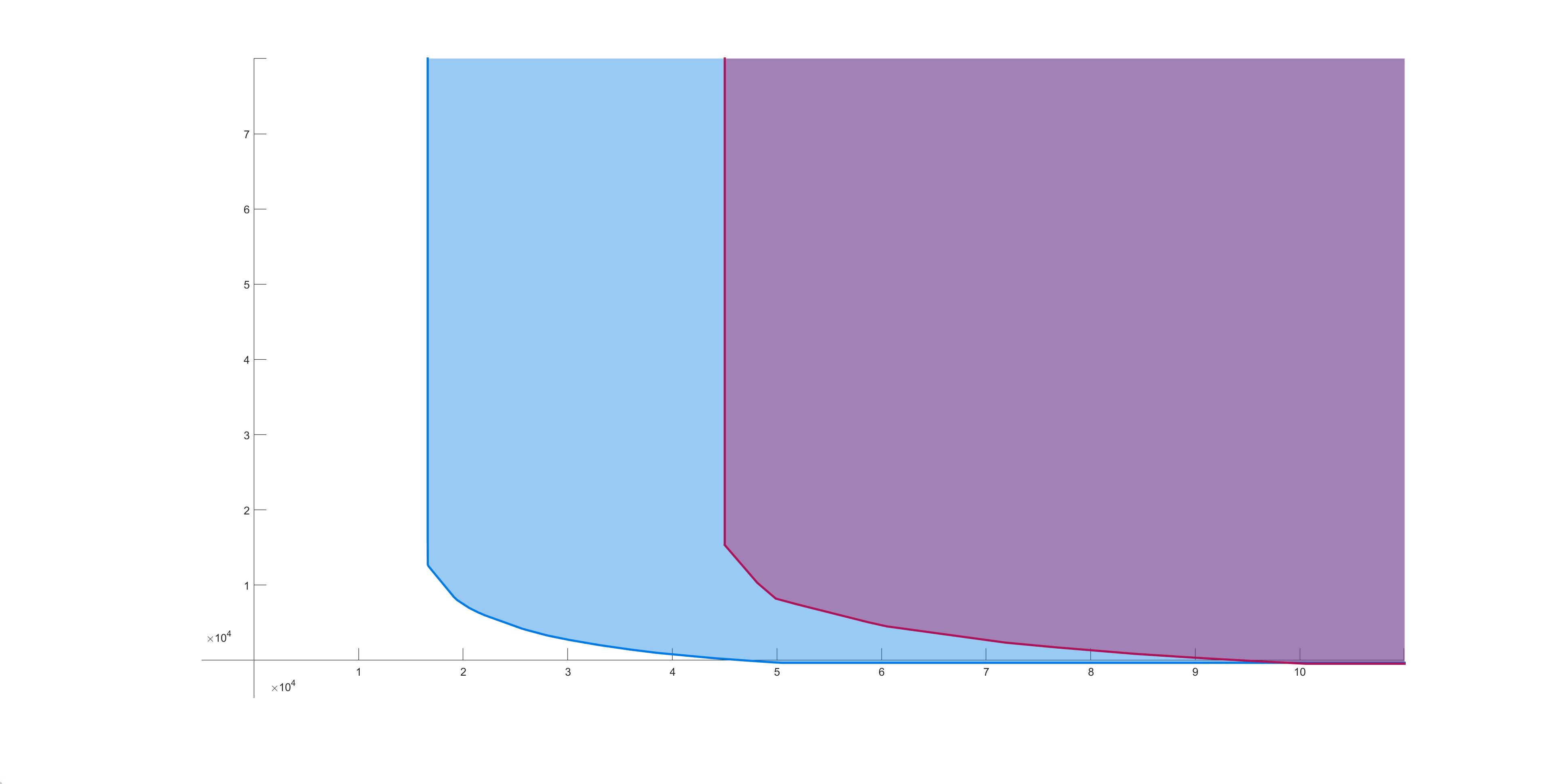}
  \caption*{$p=0.5$}
\end{minipage}
\begin{minipage}{.31\textwidth}
  \centering
  \includegraphics[width=1\textwidth]{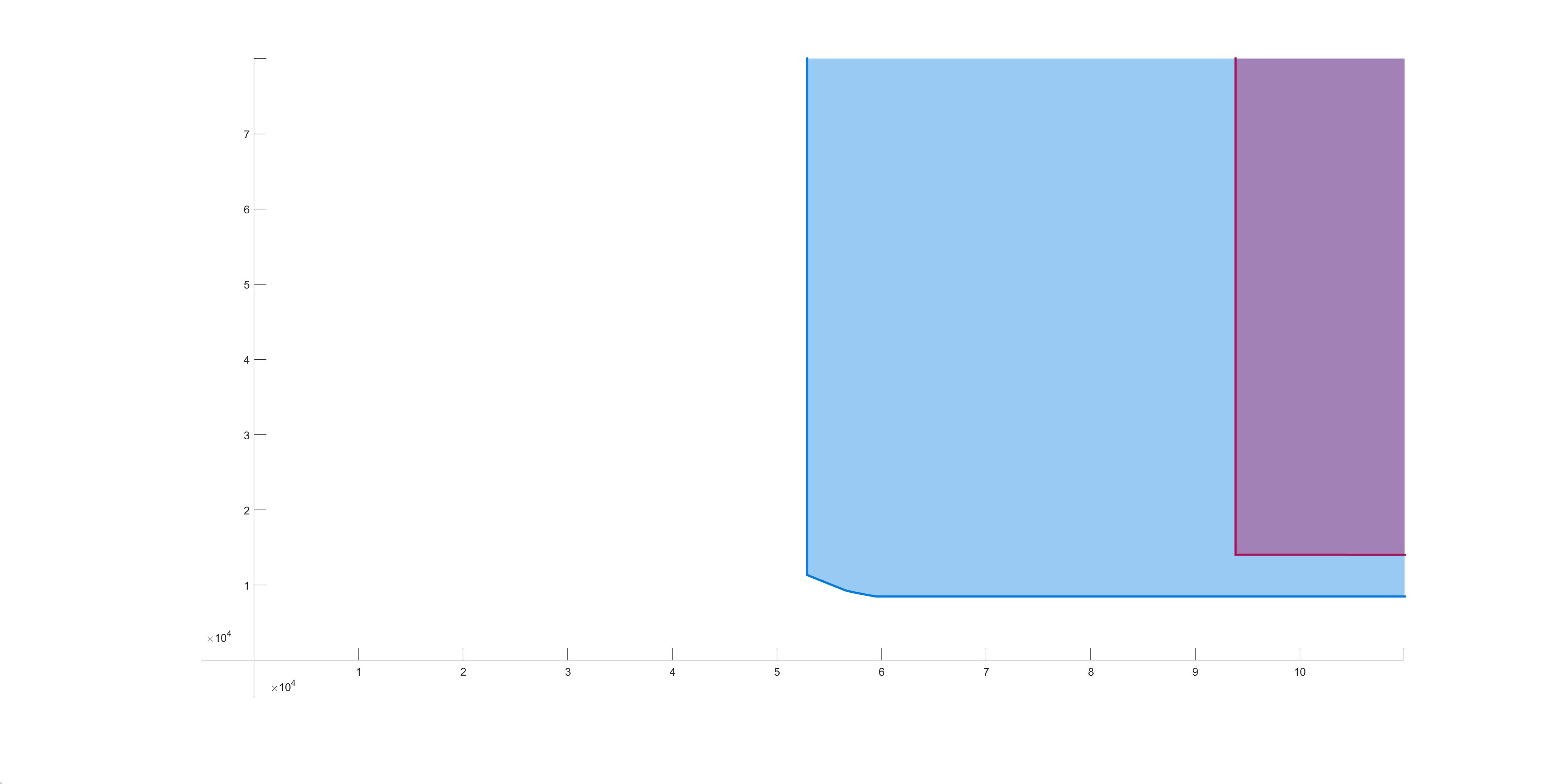}
  \caption*{$p=0.75$}
\end{minipage}
\end{figure}
The last figure shows the lower $C$-quantiles at $p=0.25$, $p=0.5$ and $p=0.75$ for the farmers without hail nets in blue and for the farmers with hail nets in red. There are two main insights that can be derived from this analysis. First, the red quantile is shifted to the right with respect to the blue quantile for all $p$. This means that the farmers with hail nets have higher yearly returns. Moreover, as $p$ increases, the horizontal distance between the quantiles increases. This implies that the combination of hail nets and insurance is more profitable the bigger the farmer. Second, the blue and red quantiles for $p=0.25$ and $p=0.5$ are aligned a little below the $x$-axes. This can be interpreted as the insurance contracts being profitable for the insurance company, but also being not too costly for the farmer. For $p=0.75$, the quantiles are substantially above the $x$-axes, indicating that bigger farmers actually make a profit out of the insurance contracts. This phenomena even increases for the farmers with hail nets.

In general, the insurance is paying off for all farmers. However, the farmers that install a hail net and stipulate an insurance are more profitable and from a certain size they have even better payoffs on the insurance.
\end{example}

\begin{example} 
\label{ExStudents}
A problem in human resource management very often is that applicants for a job (or persons potentially assigned to carry out a certain task or job) have to be evaluated according to several, often contradictory or competing criteria. The authors of this paper, as many other professionals working in academia, were subject to such evaluations. Typically, while the decision maker almost always has a clear understanding of what is better for each criterion (more publications, more project money raised, more contributions to administrative tasks etc.), the final decision is often the result of an ad hoc aggregation procedure (distribute some points for each achievement in each activity and sum the points at the end) which does not take into account the incomparableness of the different criteria and hence the different candidates. Such an aggregation is usually a scalarization: assign numbers to each candidate and then use the total order in $\R$.

Here, we use data from student results for a simple illustration what can be done in such cases and what kind of information one could expect if one would apply the methods suggested in this note. The two criteria are the average grade and the total number of credit points where the latter may serve as a proxy for the study time: the higher this number, the more courses the student finished within the time interval considered. The traditional aggregation procedure is: used the average grade as ranking for all students who earned a minimum number of credit points (made the threshold) and do not consider the others.

Note that these rankings have serious consequences: they are used to admit (or not admit) students to Erasmus programs, award prizes etc.

The average marks above 28 are highlighted in blue (one in red) in the above table. Their distribution in the right column shows that the ranking can change drastically if the cone location depth for two criteria is used instead of an (more or less arbitrarily chosen) aggregation procedure. The part highlighted in red shows that students with the same value of the cone location depth can have very different average grades. We do not claim here that the ranking according to the cone location depth is the "true" ranking (it depends on the cone which also is a choice), but we would like to point out that "traditional" decision makers should be equally aware that each ranking based on an aggregation procedure could be highly questionable and very far from being "objective."
\end{example}

\begin{figure}[H]
	\centering
	\includegraphics[width=1.1\textwidth]{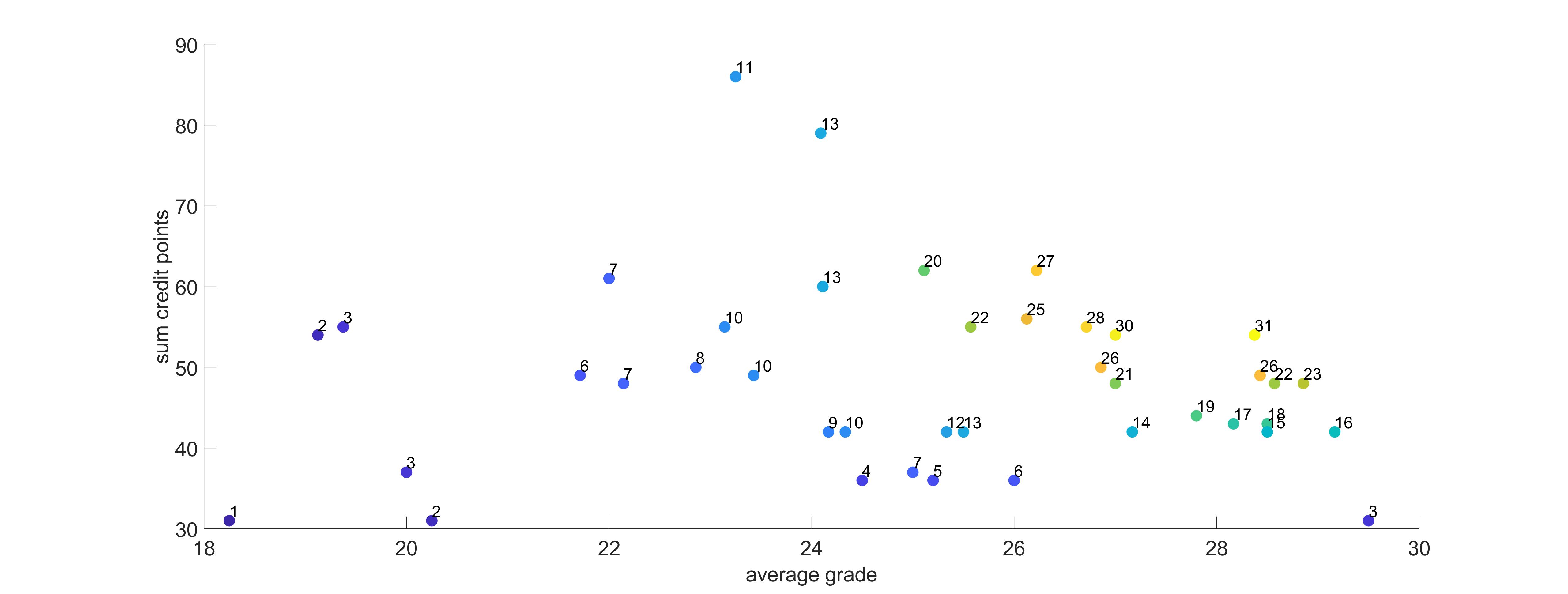}
	\caption{Each point represents a student via its average grade and the total credit points earned. The ordering cone is the $\R^2_+$. The higher the values of the cone distribution the ``better'' the student.}
	\label{fig:StudentsTourism3y}
\end{figure}

\begin{table}[htbp]
  \centering
  \begin{tabular}{@{} cccccccccccc @{}}
    \hline
    & St01 & St02 & St03 & St04 & St05 & St06 & St07 & St08 & St09 & St10 & St11 \\ 
    \hline
    CD & 31 & 30 & 28 & 27 & 26 & 26 & 25 & 23 & 22 & 22 & 21 \\ 
    AG & {\f 28.38} & 27.00 & 26.71 & 26.22 & {\f 28.43} & 26.86 & 26.13 & {\f 28.86} & {\f 28.57} & 25.57 & 27.00 \\[.25cm]
	\hline 
    & St12 & St13 & St14 & St15 & St16 & St17 & St18 & St19 & St20 & St21 & St22  \\
	\hline 
	CD & 20 & 19 & 18 & 17 & 16 & 15 & 14 & 13 & 13 & 13 & 12 \\ 
    AG & 25.11 & 27.80 & {\f 28.50} & {\f 28.17} & {\f 29.17} & {\f 28.50} & 27.17 & 25.50 & 24.11 & 24.09 & 25.33 \\[.25cm]
	\hline 
    & St23 & St24 & St25 & St26 & St27 & St28 & St29 & St30 & St31 & St32 & St33  \\
	\hline 
	CD & 11 & 10 & 10 & 10 & 9 & 8 & 7 & 7 & 7 & 7 & 6 \\ 
    AG & 23.25 & 24.33  & 23.43 & 23.14 & 24.17 & 22.86 & 25.00 & 25.00 & 22.14 & 22.00 & 26.00 \\[.25cm]
	\hline 
    & St34 & St35 & St36 & St37 & St38 & St39 & St40 & St41 & St42 &  &   \\
	\hline 
	CD & 6 & 5 & 4 & 3 & 3 & 3 & 2 & 2 & 1 &  &  \\ 
    AG & 21.71 & 25.20  & 24.50 & {\ff 29.50} & {\ff 20.00} & {\ff 19.38} & 20.25 & 19.13 & 18.25 &  &  \\[.25cm]
	\hline 
  \end{tabular}
  \caption{Student data: cone depth vs. average grades}
  \label{TabStudents}
\end{table}

\section{Conclusions}

In multivariate data analysis models, an order relation for the data points is very often (intuitively) present, but not part of the statistical analysis. In this paper, it is shown how decision makers can analyze bivariate data based on lower cone distribution functions and set-valued quantiles. It seems to us that a mere depth function approach is only appropriate if there is no (intuitive or explicitly modeled) preference present for data points they are roughly rotation symmetric. This could be subject to a debate to which this paper aims to contribute.

The algorithms for computing the new objects separate the combinatorial/nonlinear part (permutations) and linear/computational geometry part (linear rotation step for computing convex polyhedrons, see also, e.g., \cite{RousseeuwHubert17Chap} with emphasis on the link between computational geometry and depth functions/regions). This makes the computations tractable.

The approach works independently of dependence structures: in Example \ref{ExHurricanes}, the maximum windspeed and the minimal pressure are clearly strongly correlated which is not the case in Example \ref{ExStudents}. The relationships of our approach with dependence structure approaches is an interesting research question. The major difference (and difficulty) is that the lower cone distribution function is different from the joint distribution function (even if the cone is $\R^2_+$, see \cite{HamelKostner18JMVA}) which is the basic object, e.g., for copula approaches.

\medskip {\bf Funding.} D. Kostner's work was part of the project "Re-insurance of hail risks in South Tyrol" (Hail-Risk) with PI Prof. A. Wei{\ss}ensteiner funded by Free University of Bozen, Italy.

\end{document}